
\input amstex
\input epsf
\documentstyle{amsppt}
\nologo
\TagsOnRight
\def\epsfsize#1#2{\hsize}
\magnification 1200
\voffset -1cm

\def\fig#1{\bigskip\centerline{\epsffile{FIG1-#1.EPS}}\nopagebreak\medskip
\centerline{\bf Fig\. #1}\bigskip} \def\Cl#1{{\overline{#1}}} \let\x\times
\def\R{\Bbb R} \def\Z{\Bbb Z} \def\N{\Bbb N} \def\Q{\Bbb Q} \let\but\setminus
\def\section#1{\bigskip\leftline{\bf #1}\nopagebreak\medskip} \let\tl\tilde
\def\ro#1{\expandafter\def\csname#1\endcsname{{\operatorname{#1}}}}
\ro{id} \ro{im} \ro{Int} \ro{lk} \ro{rel} \let\phi\varphi \let\eps\varepsilon
\let\normal\triangleleft \let\emb\hookrightarrow \let\surj\twoheadrightarrow
\def\cal#1{\expandafter\def\csname#1\endcsname{\Cal #1}}
\cal F \cal G \cal H \cal I \cal J \cal M \cal P \cal W \cal L \cal K
\def\?{$@!@!@!@!@!$}

\topmatter 

\title $n$-Quasi-isotopy: I. Questions of nilpotence \endtitle
\author Sergey A. Melikhov and Du\v{s}an Repov\v{s} \endauthor

\abstract
It is well-known that no knot can be cancelled in a connected sum with
another knot, whereas every link can be cancelled up to link homotopy in
a (componentwise) connected sum with another link.
In this paper we address the question whether the noncancellation property
holds for some (piecewise-linear) links up to some stronger analogue of
link homotopy, which still does not distinguish between sufficiently close
$C^0$-approximations of a topological link.
We introduce a sequence of such increasingly stronger equivalence relations
under the name of $k$-quasi-isotopy, $k\in\N$; all of them are weaker
than isotopy (in the sense of Milnor).
We prove that every link can be cancelled up to peripheral structure preserving
isomorphism of any quotient of the fundamental group, {\it functorially}
invariant under $k$-quasi-isotopy; functoriality means that the isomorphism
between the quotients for links related by any allowable crossing change fits
in the commutative diagram with the fundamental group of the complement to
the intermediate singular link.
The proof invokes Baer's Theorem on the join of subnormal locally nilpotent
subgroups.
On the other hand, the integral generalized $(\lk\ne 0)$ Sato--Levine invariant
$\tl\beta$ is invariant under $1$-quasi-isotopy, but is not determined by any
quotient of the fundamental group (endowed with the peripheral structure),
functorially invariant under $1$-quasi-isotopy --- in contrast to Waldhausen's
Theorem.

As a byproduct, we use $\tl\beta$ to determine the image of the Kirk--Koschorke
invariant $\tl\sigma$ of fibered link maps.
\endabstract

\address Steklov Mathematical Institute, Division of Geometry and Topology;
\newline ul. Gubkina 8, Moscow 119991, Russia \endaddress
\curraddr University of Florida, Department of Mathematics;\newline
358 Little Hall, PO Box 118105, Gainesville, FL 32611-8105, U.S. \endcurraddr
\email melikhov\@math.ufl.edu, melikhov\@mi.ras.ru \endemail
\address Institute of Mathematics, Physics and Mechanics, University of
Ljubljana, 19 Jadranska cesta, 1001 Ljubljana, Slovenia \endaddress
\email dusan.repovs\@fmf.uni-lj.si \endemail
\endtopmatter

\document 

\head 1. Introduction \endhead

\subhead \oldnos1.\oldnos1. The definition \endsubhead
Consider PL links $L,L'\:mS^1\emb S^3$, where $mS^1$ denotes
$S^1_1\sqcup\dots\sqcup S^1_m$, and let $f\:mS^1\to S^3$ be any singular
link in a generic PL homotopy between them.
Then $f$ is called a {\it $0$-quasi-embedding}, or simply a {\it link map},
if its (unique) double point is a self-intersection of some component
$f(S^1_i)$ (rather than an intersection between distinct components).
If this is the case, $f(S^1_i)$ is a wedge of two loops $\ell,\ell'$
(following \cite{KL1}, we call them the {\it lobes} of $f(S^1_i)$), and
we say that $f$ is a {\it $1$-quasi-embedding,} if at least one of them, say
$\ell$, is null-homotopic in the complement to the non-singular components,
i.e\. in $X:=S^3\but f(mS^1\but S^1_i)$.
If this also is the case, we say that $f$ is a {\it $2$-quasi-embedding,} if
a generic PL null-homotopy $F\:D^2\to X$ of $\ell=F(\partial D^2)$ can be
chosen so that for some arc $I\i\ell'$ containing $\ell'\cap F(D^2)$, the
polyhedron $F(D^2)\cup I$ is null-homotopic in $X$.

In general, a PL map $f\:mS^1\to S^3$ with precisely one double point $x$
will be called a {\it $k$-quasi-embedding}, $k\in\N$, if, in addition to
the singleton $P_0=\{x\}$, there exist compact subpolyhedra
$P_1,\dots,P_k\i S^3$ and arcs $J_0,\dots,J_k\i mS^1$ such that
$f^{-1}(P_j)\i J_j$ for each $j\le k$, $P_j\cup f(J_j)\i P_{j+1}$ for each
$j<k$, and the latter inclusion is null-homotopic for each $j<k$.
As a matter of convenience, we shall also assume that
$\partial J_0=f^{-1}(x)$.

A reader familiar with embeddings of higher dimensional manifolds may
recognize this construction as $k$ steps of the
Penrose--Whitehead--Zeeman--Irwin trick \cite{PWZ}, \cite{Ir}, \cite{Ze1},
adjusted to codimension two; see also \cite{St2; Comment~ 6.3}.
Note that the original construction of Casson handles \cite{Ca2} can be
thought of as an instance of infinitely many steps of this (adjusted) trick.

Two links $L,L'\:mS^1\to S^3$ are said to be {\it $k$-quasi-isotopic,} if
they can be joined by a generic PL homotopy, all whose singular links are
$k$-quasi-embeddings.
Thus $0$-quasi-isotopy is nothing but {\it link homotopy}
(i.e\. homotopy through link maps), whereas {\it PL isotopy} (i.e\. PL homotopy
through one-to-one maps) implies $k$-quasi-isotopy for all $k\in\N$.
(Note that a PL isotopy may insert/shrink local knots.)
It is not hard to see (cf\. {\it ``$n$-Quasi-isotopy II''}) that
the $(k+1)^{\text{th}}$ Milnor link and the $(k+1)^{\text{th}}$
Whitehead link are $k$-quasi-isotopic to the unlink.

It will be convenient to have two modifications of the definitions.
{\it Strong (weak) $k$-quasi-embedding}, and consequently {\it strong (weak)
$k$-quasi-isotopy}, are defined as above, except that ``compact polyhedra'' is
replaced with ``PL $3$-balls'' (resp\. ``null-homotopic'' is replaced with
``induces trivial homomorphisms on reduced integral homology'').
Strong $k$-quasi-isotopy may be thought of as arising from a slightly
different interpretation of the PWZI trick, where the engulfing lemma is
taken for granted, rather than regarded as a part of the construction.

\subhead \oldnos1.\oldnos2. A motivation \endsubhead
The concept of $k$-quasi-isotopy arose in attempts of the first author to
solve the following problem (cf\. \cite{Me1; Questions I and III}).
Let $M^1$ be a compact $1$-manifold (it does not seem to really matter which
one).
A (continuous) map $f\:M\to\R^3$ is said to be {\it realizable by TOP
isotopy}, if there exists a homotopy $h_t\:M\to\R^3$, $t\in [0,1]$, such that
$h_1=f$ and each $h_{t_0}$ with $t_0<1$ is one-to-one.

\proclaim{Problem 1.1} Does there exist a map $M^1\to\R^3$, non-realizable
by TOP isotopy?
\endproclaim

As pointed out in \cite{Me1}, it would be especially plausible to find such
a map among locally flat TOP immersions (i.e\. maps that tamely embed
a neighborhood of every point).
For each $n\ge 3$ there exists a locally flat TOP immersion $S^n\to\R^{2n}$,
which is $\eps$-approximable by an embedding for each $\eps>0$, but is not
realizable by TOP isotopy \cite{Me3} (see also \cite{Me1; Theorem 1.12}).

\fig 1

A possible candidate $f\:I\sqcup I\to\R^3$ (with one double point)
for being non-realizable by TOP isotopy is depicted in Fig\. 1.
This locally flat TOP immersion can be regarded as a ``connected sum'' of
infinitely many (left handed) Whitehead links (resembling the Wilder arc,
which is an infinite ``connected sum'' of trefoils).
Note also that $f$ is a composition of a PL map and a TOP embedding; it
was shown in \cite{Me1} that in codimension $\ge 3$, such a composition is
realizable by a TOP isotopy if and only if it is $\eps$-approximable by
an embedding for each $\eps>0$.

\proclaim{Proposition 1.2} The map $f$ of Fig\. 1 is not realizable by
PL isotopy.
\endproclaim

The proof is not hard, but reveals the link-theoretic nature of Problem 1.1.
We first introduce some notation to be used in this proof and throughout \S2.
\smallskip

A two-component {\it string link} is a PL embedding
$L\:(I_+\sqcup I_-,\partial)\emb (I\x\R^2,\partial)$ such that
$L(i,\pm)=(i,\pm p)$ for $i=0,1$ and some fixed $p\in\R^2\but\{0\}$.
String links are considered up to ambient isotopy, and their connected sum
$L_1\#L_2$ is defined by stacking together the two copies of $(I_+,I_-,I)$.

If $L$ is a two-component string link, let $denom(L)\:S^1\emb S^3$ be the
``denominator'' closure of $L$, obtained by adjoining to $L$ the two arcs
$\partial I\x\{\lambda p\mid\lambda\in[-1,1]\}$ in $\partial I\x\R^2$,
and let $numer(L)\:S^1\sqcup S^1\emb S^3$ be the ``numerator'' closure of
$L$, defined to be the composition of
$L/\partial\:I_+/\partial\sqcup I_-/\partial\emb I/\partial\x\R^2$
and the inclusion $S^1\x\R^2\i S^1\x D^2\cup D^2\x S^1=\partial(D^2\x D^2)$.

Finally, we call the string link $W$ that appears three times in Fig\. 1
the {\it Whitehead string link}, since its linking number zero closure
(cf\. \S2) is the Whitehead link $\W$.

\demo{Proof}
Suppose that $h_t\:I\sqcup I\to\R^3$, $t\in [0,1]$, is a homotopy such that
$h_1=f$, each $h_{t_0}$ with $t_0<1$ is one-to-one, and which restricts to
a PL homotopy for $t\in [0,1)$.
We may assume that $f(I\sqcup I)$ is contained in the interior of
$I\x\R^2\i\R^3$.
Let $U$ be a regular neighborhood of the $4$ endpoints
$f(\partial I\sqcup\partial I)$ in $I\x\R^2$ such that $f^{-1}(U)$
is a regular neighborhood of $\partial I\sqcup\partial I$ in $I\sqcup I$,
and each of the $4$ balls of $U$ meets $\partial I\x\R^2$ in a disk.
Without loss of generality, $h_{t_0}(I\sqcup I)$ is contained in the interior
of $I\x\R^2$ and $h_{t_0}(\partial I\sqcup\partial I)$ is contained in $U$
for each $t_0$.
Then $h_t$ for $t\in [0,1)$ extends, via an identification of $I$ with
a subarc of $I_\pm$, to a PL isotopy
$H_t\:(I_+\sqcup I_-,\partial)\emb (I\x\R^2,\partial)$, keeping the boundary
fixed and such that $H_t(I_+\but I\sqcup I_-\but I)\i U$.

Each $H_{1-\eps}$ is a string link, which for a sufficiently small $\eps>0$
can be represented as a connected sum $W\#\dots\# W\# L'$ of arbitrarily many
copies of the Whitehead string link $W$ and some residual string link $L'$.
Therefore, in order to get a contradiction it suffices to find an invariant
$\I$ of PL isotopy of $2$-component string links such that

(i) $\I(L)$ is a nonnegative integer for every $L$;

(ii) $\I(W)>0$;

(iii) $\I(L_1\# L_2)\ge\I(L_1)+\I(L_2)$ for any $L_1$ and $L_2$.

\noindent
If $L$ is a two-component string link, let $loc(L)$ denote the knot product
of the local knots of $L$ or, more precisely, of all one-component factors
in the prime factorization \cite{Has} of the non-split $3$-component link
$numer(L)\cup M$, where $M$ is a meridian of the solid torus $S^1\x D^2$
from the definition of $numer(L)$ above.
A collection of disjoint balls $B_i$ in $S^3\but M$ meeting $numer(L)$ in
arcs $I_i$ such that $(B_i,I_i)$ represent the one-component prime factors of
$numer(L)\cup M$ can be chosen so that each $B_i$ is disjoint from a given
$2$-disk $D$ bounded by $M$ and meeting $numer(L)$ in two points (this is
proved similarly to the uniqueness of knot factorization, see \cite{Lic}).
It follows that $loc(L)$ is a factor of $denom(L)$, i.e\. there exists a knot
$glob(L)$ such that $denom(L)=loc(L)\# glob(L)$.

We define $\I(L)$ to be the genus of $glob(L)$.
Then $\I(L)$ is invariant under PL isotopy since $glob(L)$ is.
Since $loc(W)$ is the unknot and $denom(W)$ is the trefoil, $\I(W)=1$.
Finally, $loc(L_1\#L_2)=loc(L_1)\#loc(L_2)$ for any $L_1$, $L_2$ by a version
of the above assertion where the $2$-disk $D$ is replaced by two such $2$-disks
with disjoint interiors (the proof is analogous).
Therefore, since knot genus is additive under connected sum,
$\I(L_1\#L_2)=\I(L_1)+\I(L_2)$. \qed
\enddemo

\remark{Remark} One could try to prove Proposition 1.1 only using closed links.
As in the above proof, it is easy to extend $h_t$ for $t\in [0,1)$ to
a PL isotopy $H_t\:S^1\sqcup S^1\emb\R^3$ such that the link $H_{1-\eps}$
for a sufficiently small $\eps>0$ is ambient isotopic a connected sum%
\footnote{Following \cite{Co2; \S5} and \cite{Hi; \S2.7}, we denote by
$L_1\# L_2$ any link of the form $numer(S_1\#S_2)$, where $L_1=numer(S_1)$
and $L_2=numer(S_2)$.
This multi-valued operation of componentwise connected sum should not be
confused with the multi-valued operation of connected sum along selected
components as in the prime decomposition theorem of \cite{Has}.}
$\W\#\dots\#\W\#\Cal L$ of arbitrarily many copies of the Whitehead link
$\W$ and some residual link $\Cal L$.
However, it turns out that there exists no invariant of $2$-component links
satisfying the analogues of conditions (i)--(iii); indeed, a certain
connected sum of the Whitehead link with its reflection (along the bands,
``orthogonal'' to the ``mirror'') is ambient isotopic to the unlink.
\endremark
\bigskip

Problem 1.1 would be solved if we could show that the invariant $\I$ from
the proof of Proposition 1.2 stabilizes for sufficiently close (in the $C^0$
topology) embedded PL approximations of any given topological string link;
for this would imply invariance of $\I$ under TOP isotopy with parameter in
the (compact) interval $[0,1]$.
The following example shows that this is not the case, moreover, $\I$ assumes
arbitrarily large values on arbitrarily close PL approximations of a certain
wild string link.

\fig 2

\remark{Example} Consider the wild link $\L(2S^1)=C\cup\K$, depicted in
Fig\. 2 (ignore the dashed lines and the disk for now), and let $K$ be
a PL knot that is the image of an $\eps$-approximation of the parametrization
of $\K$ for some small $\eps>0$.
As a knot, $K$ can be gradually untangled from left to right into a connected
sum of trefoils $(B^3_i,B^3_i\cap K)$, $1\le i\le N$, and possibly some
residual prime knots $(B^3_i,B^3_i\cap K)$, $N+1\le i\le M$, where
$N=N(\eps)\to\infty$ as $\eps\to 0$, and the balls $B^3_i$ are pairwise
disjoint.
None of these $N$ trefoils is a local knot (i.e\. a one-component factor in
the prime decomposition) of the link $L:=C\cup K$, since the free homotopy
class of $C$, regarded as a conjugacy class in the amalgamated free product
$\pi_1(S^3\but K)\cong
\pi_1(B^3_1\but K)\underset\Z\to*\cdots\underset\Z\to*\pi_1(B^3_M\but K)$,
can be seen to remain nontrivial under each homomorphism
$\pi_1(S^3\but K)\to\pi_1(B_i^3\but K)$, $i\le N$, given by abelianization
of all the factors except for the $i^{\text{th}}$ one.%
\footnote{Thus $\L$ is a counterexample to the assertion in the last sentence
of Example 1.3 in \cite{Me1}.
The mistake in the proof of this assertion was pointed out in \cite{Me2}.}
This homomorphism does not depend on the choice of the prime decomposition
$(B^3_1,\dots,B^3_M)$ of $K$ since the proof of the uniqueness of knot
factorization (see \cite{BZ}) shows that the collection of knot exteriors
$(\Cl{B^3_1\but N(K)},\dots,\Cl{B^3_M\but N(K)}\,)$, where $N(K)$ is a small
regular neighborhood of $K$, is unique up to a self-homeomorphism of
$(S^3,K)$.

Now let $\K'$ and $K'$ be obtained from $\K$ and $K$ by the clockwise twist
by $2\pi$ of the visible side of the disk shown in Fig\. 2 (the disk is
bounded by $C$ and meets $K$ in two points); set $\L'=C\cup\K'$ and
$L'=C\cup K'$.
The twist yields a homeomorphism between the complements of $L$ and $L'$
sending meridians of $K$ onto those of $K'$, therefore $loc(L')=loc(L)$.
If $L'=numer(S)$ for some string link $S$, the knot $denom(S)$ is a band
connected sum of the components of $L'$; moreover, $S$ is uniquely determined
by $L'$ and the band.
The band indicated by dashed arcs in Fig\. 2 yields $denom(S)=K$.
Since $loc(S)$ is a factor of $loc(L')$, the knot $glob(S)$ has at least $N$
trefoils among its prime factors, and so $\I(S)\ge N$.
Thus $\I$ can be arbitrarily large for arbitrarily close PL approximations of
the wild string link $\Cal S$ that is determined by $\L'$ and the band.
\endremark

\subhead \oldnos1.\oldnos3. The main problem \endsubhead
In general, it seems to be not easy to find a nonzero invariant of topological
string links, satisfying the analogues of conditions (i) and (iii) from
the proof of Proposition 1.2 (compare \cite{Me1; Example 1.5} and \cite{KY}).
Since we are not interested in exploring invariants of wild links per se,
we will take a different approach in this paper, staying in the class of
PL links by the price of imposing an additional restriction on the desired
invariant, which ensures its stabilization on sufficiently close approximations
of any topological string link.
This approach is based on

\proclaim{Theorem 1.3} For each $k\in\N$ and any TOP isotopy
(i.e\. homotopy through one-to-one maps) $h_t\:mS^1\emb\R^3$, $t\in [0,1]$,
there exists an $\eps>0$ such that any generic PL homotopy, $\eps$-close to
$h_t$, is a strong $k$-quasi-isotopy.
\endproclaim

\demo{Proof} For the sake of clarity we only give an explicit proof for $k=1$.

Let $\delta>0$ be the minimal distance between $h_t(S^1_i)$ and $h_t(S^1_j)$,
$i\ne j$, $t\in I$.
Let $C_i\i (S^1_i)^2\x I$ denote the set of all triples $(p,q,t)$ such that
$p$ and $q$ are not antipodal in $S^1_i$, and the $h_t$-image of the shortest
arc $[p,q]\i S^1_i$ lies in the open $\frac{\delta}5$-neighborhood of
the middle point $m_{pq}$ of the line segment $[h_t(p),h_t(q)]\i\R^3$.
Then $C_i$ contains a neighborhood of $\Delta_{S^1_i}\x I$, and it follows
that the union $C$ of all $C_i$'s contains $\cup_t (h_t^2)^{-1}(N)$ for some
neighborhood $N$ of $\Delta_{\R^3}$.
Hence the minimal distance $\gamma$ between $h_t(p)$ and $h_t(q)$ over all
triples $(p,q,t)$ in the (compact) complement to $C$ is non-zero.
Set $\eps=\min(\frac{\gamma}3,\frac{\delta}5)$.

Let $l_t$ be a generic PL homotopy, $\eps$-close to $h_t$.
Then $l_t$ has at most one double point $x:=l_t(p)=l_t(q)$ for each $t\in I$.
In this case $h_t(p)$ and $h_t(q)$ are $\frac{2\gamma}3$-close.
So $(p,q,t)$ lies in $C_i$ for some $i$.
Hence $h_t([p,q])$ lies in the $\frac{\delta}5$-neighborhood of $m_{pq}$.
It follows that the loop $\ell:=l_t([p,q])$ lies in the
$\frac{3}5\delta$-neighborhood of $x$.
We claim that $\ell$ is null-homotopic in $S^3\but l_t(mS^1\but S^1_i)$.
Indeed, if the radial null-homotopy of $\ell$ towards $x$ met $l_t(S^1_j)$,
then a point of $\ell$ and a point of $l_t(S^1_j)$ would lie on the same
radius of the $\frac{3}5\delta$-neighborhood of $x$, consequently a point
of $h_t([p,q])$ and a point of $h_t(S^1_j)$ would be $\delta$-close,
contradicting our choice of $\delta$. \qed
\enddemo

The same argument works to prove the analogous assertion for string links.
(The definition of $k$-quasi-isotopy for string links is the same, except
that ``at least one of the two loops'' must be changed to ``the loop'' in
the reformulations for $k=1,2$.)

\proclaim{Corollary 1.4}
(a) For each $k\in\N$, strong $k$-quasi-isotopy classes of all sufficiently
close PL approximations to a topological link (or topological string link)
coincide.

(b) TOP isotopic PL links (or PL string links) are strongly
$k$-quasi-isotopic for all $k\in\N$.
\endproclaim

Thus, positive solution to Problem 1.1 is implied by positive answer to
either part of the following question of independent interest.

\proclaim{Problem 1.5} Does there exist, for some $k,m\in\N$, a nonnegative
integer valued invariant $\I$ of $k$-quasi-isotopy of $m$-component PL string
links such that

(a) $\I(L\# L')\ge\I(L)+\I(L')$ for any $L$, $L'$, and $\I(L_0)>0$ for
some $L_0$, preferably with unknotted components?

(b) {\rm (non-commutative version)} $\I(L\# L')\ge\I(L)$ for any $L$, $L'$,
and there exist $L_0,L_1,\dots$, preferably with unknotted components,
such that $\I(L_0\#\dots\#L_n)\ge n$ for each $n$?

(a$'$), (b$'$) Same questions for PL links (with any choice of bands in
$L\# L'$ and some choice of bands in $L_0\#\dots\#L_n$).
\endproclaim

Additional interest in Problem 1.5 arises from its interpretation as a proper
formalization of the intuitive question, whether the ``phenomenon of linking
(modulo knotting)'', regarded solely as a feature of interaction between
distinct components, admits the possibility of ``accumulation of complexity''
(like it happens in the case of the ``phenomenon of knotting'').

For $k=0$ the answer to Problem 1.5 is well-known to be negative since
concordance implies link homotopy \cite{Gi}, \cite{Go}.
Indeed, the connected sum of any PL link $L$ with its reflection $\rho L$
(i.e\. the composition of $L$ with reflections in both domain and range)
along the bands, ``orthogonal'' to the ``mirror'', is a slice link.
For completeness, the proof of the analogous assertion for string links is
included in the proof of Theorem 2.2 below.

\subhead \oldnos1.\oldnos4. Main results \endsubhead
A straightforward way to find an invariant as required in Problem 1.5 is to
extract it from some quotient of the fundamental group
$\pi(L):=\pi_1(S^3\but L(mS^1))$, invariant under $k$-quasi-isotopy.
In view of Stallings' Theorem on lower central series \cite{St1},
\cite{Ca1}, \cite{Co1} this is impossible if the quotient is nilpotent.
Indeed, it follows from the Stallings Theorem that any nilpotent quotient is
invariant under link concordance, whereas any PL link (or string link) is
cancelled by its reflection up to concordance.
This suggests that the first step towards the above problems might be to look
for non-nilpotent quotients of $\pi(L)$, invariant under $k$-quasi-isotopy of
$L$.

The following encouraging fact exploits the power of the generalized
Sato--Levine invariant $\tl\beta$, whose remarkability emerged independently
in the works of Traldi \cite{Tr; \S10}; Polyak and Viro \cite{PV} (see also
\cite{\"Os}, \cite{AMR}); Kirk and Livingston \cite{KL1}, \cite{Liv};
Akhmetiev \cite{Ah}; \cite{AR}; Nakanishi and Ohyama \cite{NO}; and,
implicitly (see proof of Theorem 2.6 below), Koschorke
\cite{Ko; discussion of Fact 2.9}.

\proclaim{Corollary 2.3} No nilpotent quotient of the fundamental group
(even if endowed with the peripheral structure%
\footnote{We recall that the {\it peripheral structure} in the fundamental
group of an $m$-component link $L$ is a collection of $m$ pairs $(m_i,l_i)$,
defined up to simultaneous conjugation of $m_i$ and $l_i$ by an element of
$\pi(L)$, where $m_i$ is a meridian to the $i^{\text{th}}$ component of $L$
and $l_i$ is the corresponding longitude.
The {\it meridian} and the {\it longitude} (or parallel) corresponding to
a path $p$ that joins the basepoints of the link complement and of
the boundary of a small regular neighborhood $T_i$ of the $i^{\text{th}}$
component, are the images of the distinguished generators of
$\pi_1(\partial T_i)\simeq H_1(\partial T_i)\simeq H_1(\Cl{S^3\but T_i})
\oplus H_1(T_i)$ under the homomorphism $i_*\:\pi_1(\partial T_i)\to\pi(L)$
determined by $p$.}%
) can be a complete invariant of $1$-quasi-isotopy of $2$-component links.

Similarly for string links; specifically, the connected sum $W\# rW$ of
the Whitehead string link and its reflection is not $1$-quasi-isotopic to
the trivial string link.
\endproclaim

However, we were unable to find any non-nilpotent quotient of the fundamental
group, invariant under $k$-quasi-isotopy for some $k$.
The difficulty that we encountered can be summarized as follows.

Let $R$ be an equivalence relation on links, obtained from ambient isotopy
by allowing certain types of transversal self-intersections of components,
and $\F$ be a function assigning to every PL link $L$ some group $\F(L)$
together with an epimorphism $p_L\:\pi(L)\to\F(L)$.
Then we say that $\F$ is {\it functorially invariant} under $R$ with
respect to $p_L$ if it is invariant under ambient isotopy, and for any links
$L_+$, $L_-$ obtained one from another by an allowed crossing change, with
$L_s$ being the intermediate singular link, there is an isomorphism between
the groups $\F(L_+)$ and $\F(L_-)$ making the following diagram
commutative:
$$\CD
\pi_1(S^3\but\text{regular neigh}\hskip -20pt@.\hskip-20pt
\text{borhood of }L_s(mS^1))\\
@Vi_*VV@VVi_*V\\
\pi(L_+)@.\pi(L_-)\\
@Vp_{L_+}VV@VVp_{L_-}V\\
\F(L_+)@>\simeq>>\F(L_-)
\endCD$$

\proclaim{Corollary 3.7} Every quotient of the fundamental group of a link
(or string link), functorially invariant under $k$-quasi-isotopy for some
$k$, is nilpotent.
\endproclaim

By a theorem of Waldhausen (see \cite{Ka}) the ambient isotopy class of $L$ is
completely determined by the fundamental group $\pi(L)$ together with
the peripheral structure $P(L)$.
Therefore, given an equivalence relation $R$ on links, it must be possible, in
principle, to extract any invariant of links up to $R$ from the pair $(\pi,P)$.
Nevertheless, the two results above combine to produce

\proclaim{Corollary 1.6} Let $\G_k(L)$ denote the finest quotient of $\pi(L)$,
functorially invariant under $k$-quasi-isotopy%
\footnote{To see that $\G_k(L)$ is well-defined, note that it is the quotient
of $\pi(L)$ over the product of the kernels of
$i_*\:\pi(L)\to\pi_1(S^3\x I\but H(mS^1\x I))$ for all $k$-quasi-isotopies $H$
starting with $L$.}%
, and $\P_k(L)$ the finest peripheral structure in $\G_k(L)$, invariant under
$k$-quasi-isotopy%
\footnote{More precisely, $\P_k(L)$ denotes the collection of $m$ pairs
$(\bar m_i,\Lambda_i)$, each defined up to conjugation, where
$\bar m_i\in\G_k(L)$ is the coset of a meridian $m_i\in\pi(L)$ to
the $i^{\text{th}}$ component of $L$, and $\Lambda_i\i\G_k(L)$ is the set of
cosets $\bar l_{i\alpha}\in\G_k(L_\alpha)\simeq\G_k(L)$ of the longitudes
$l_{i\alpha}$ corresponding to some representatives $m_{i\alpha}$ of $\bar m_i$
in the fundamental groups of all links $L_\alpha$, $k$-quasi-isotopic to $L$
(that $m_{i\alpha}$ are meridians follows from functoriality).}%
. For each $k>0$, the pair $(\G_k,\P_k)$ is not a complete invariant of
$k$-quasi-isotopy, even for $2$-component links (and string links).
\endproclaim

In fact, the gap between $k$-quasi-isotopy and what this pair can tell about
it must be extremely wide, apparently including distinction of most slice
and boundary links if $k$ is large enough.
Certainly, concordant links need not be $1$-quasi-isotopic (see \S2), whereas
the boundary of two Whitehead-linked M\"obuis bands \cite{NS; Fig\. 7} is
a $\Z/2$-boundary link with $\lk=0$, which is $1$-quasi-isotopically
non-trivial since it has a nonzero Sato--Levine invariant (cf\. \S2).

\subhead \oldnos1.\oldnos5. Discussion and remarks \endsubhead
Surprisingly, for $k=0$ it seems to be a hard problem whether $(\G_0,\P_0)$ is
a complete invariant of link homotopy (= $0$-quasi-isotopy).
It is not hard to see (see \S3 or \cite{Le; proof of Theorem 4}) that
$\G_0(L)$ coincides with Milnor's link group $\G(L)$ (i.e\. the quotient of
$\pi(L)$ obtained by forcing each meridian to commute with all of its
conjugates \cite{Mi}), and $\P_0(L)$ with the peripheral structure $\P(L)$
where the set $\Lambda_i$ corresponding to the class $\bar m_i$ of a meridian
$m_i\in\pi(L)$ is the right coset $N(m_i)\bar l_i$ of the subgroup
\footnote{We use the notation $[a,b]=a^{-1}b^{-1}ab$, $a^b=b^{-1}ab$ throughout
the paper, but in the notation $[a,b]=aba^{-1}b^{-1}$, $a^b=bab^{-1}$ of
\cite{Le} the same formula happens to define the same group (by substituting
$g^{-1}$ for $g$).}
$N(m_i)=\left<[g^{-1},g^{m_i}]\mid g\in\G(L)\right>$ containing the class
$\bar l_i$ of the longitude corresponding to $m_i$.
In fact, Levine explicitly considered only a coarser peripheral structure
$\P'(L)$ in \cite{Le}, with the right coset of $N(m_i)$ being enlarged to the
coset of the normal closure $N_i$ of $N(m_i)$; but this $\P'$ was still finer
than Milnor's original $\P''$ where instead of $N_i$ the normal closure of
$m_i$ in $\G$ was used \cite{Mi}.

It has been known for some time that Milnor's pair $(\G,\P'')$ is a complete
invariant of link homotopy in the cases of homotopically Brunnian and
$\le 3$-component links \cite{Mi}, as well as of links obtained by
adjoining a component in an arbitrary way to a homotopically Brunnian link
\cite{Le}; yet another case was found in \cite{Hu2; Theorem 3.1}.
This is not the case for $4$-component links \cite{Hu1}, where nevertheless
Levine's peripheral structure $\P'$ yields a complete classification
\cite{Le}.
On the other hand, Levine's program \cite{Le} for inductive determination of
the link homotopy class of a link from its $(\G,\P')$ via a description of
the automorphisms of the Milnor group of a principal sublink, induced by its
self-link-homotopies, was found in \cite{Hu1}, \cite{Hu2} to be not as certain
to succeed as it might have been expected, because for certain $5$-component
link not all automorphisms preserving the peripheral structure $\P'$ are
induced by self-link-homotopies \cite{Hu1}.
Another viewpoint of these difficulties is that they arise from non-uniqueness
of representation of link as the closure of a string link \cite{HL}; indeed,
string links are completely classified by their (integer) $\bar\mu$-invariants
and so by the string link analogue of $(\G,\P'')$.
\smallskip

We conclude this section with general remarks on the notion of
$k$-quasi-isotopy.

\remark{Remarks. (i)}
The term ``$k$-quasi-isotopy'' is motivated by the following considerations.
A compact polyhedron $X$ is said to be {\it quasi-embeddable} (cf\.
\cite{SSS}) into a PL manifold $Q$, if for each $\eps>0$ there exists a map
$X\to Q$ with all point-inverses of diameter at most $\eps$ (such a map is
called an {\it $\eps$-map}).
It is natural to call two embeddings $X\emb Q$ {\it quasi-isotopic} if for
each $\eps>0$ they can be joined by a homotopy in the class of $\eps$-maps.
(For links of spheres in a sphere in codimension $\ge 2$ this is clearly
equivalent to being link homotopic.)
These correspond to the case $k=0$ of the following definition.
\endremark

\remark{(ii)} In order to clarify our motivations (and for future reference)
we define $k$-quasi-isotopy ($k\in\N$) in the most general situation,
that is for PL embeddings of a compact polyhedron $X$ into a polyhedron $Y$.
Let $f\:X\to Y$ be an arbitrary PL map, and consider the singular set
$S_f=\text{closure}(\{x\in X\mid f^{-1}(f(x))\ne\{x\}\})$.
Then $f$ is called a {\it strong $k$-quasi-embedding} if there exist
a commutative PL diagram

\def\s{\ \ \surj\ \ }
$$\CD
S_f       @.\i@.J_0       @.        @.\i      @.  @.J_1       @.
 @.\i      @.  @.\cdots@.  @.J_{k-1}       @.        @.\i      @.  @.J_k
  @.\i         @.X\\
@VVV          @.|         @.\searrow@.        @.  @.|         @.\searrow
 @.        @.  @.      @.  @.|             @.\searrow@.        @.  @.|
  @.           @.@VfVV\\
P_0       @.  @.\i        @.        @.P_1     @.  @.\i        @.
 @.P_2     @.  @.\cdots@.  @.\i            @.        @.P_k     @.  @.\i
  @.           @.Y\\
@VV\phi_0V@.  @VV\psi_0V  @.        @VV\phi_1V@.  @VV\psi_1V  @.
 @VV\phi_2V@.  @.      @.  @VV\psi_{k-1}V  @.        @VV\phi_kV@.  @VV\psi_kV
  @.           @VVV\\
Q_0       @.\s@.\tilde Q_0@.\s      @.Q_1     @.\s@.\tilde Q_1@.\s
 @.Q_2     @.\s@.\cdots@.\s@.\tilde Q_{k-1}@.\s      @.Q_k     @.\s@.\tilde Q_k
  @.\ \ \to\ \ @.pt\\
\endCD$$
\bigskip

\noindent
and a triangulation
$K_0\surj\tilde K_0\surj K_1\surj\tilde K_1\surj\dots\surj K_k\surj\tilde K_k$
of the bottom line, satisfying the following conditions for each $i=0,\dots,k$:

\smallskip
(a) $f^{-1}(P_i)\i J_i$;

(b) $\phi_i^{-1}(C)$ is collapsible for each cone $C$ of the dual cone
complex $K_i^*$;

(c) $\psi_i^{-1}(C)$ is collapsible for each cone $C$ of the dual cone
complex $\tilde K_i^*$.
\smallskip

\noindent
Note that $\phi_0\:P_0\to Q_0$ can always be chosen to be
$\id\:f(S_f)\to f(S_f)$; and that if $X$ and $Y$ are PL manifolds,
``is collapsible'' may be replaced by ``is a codimension zero PL ball'' by
the theory of regular neighborhoods.
(Recall that the empty set is not collapsible.)
The definition of a {\it $k$-quasi-embedding} is similar, with (b) and (c)
replaced by the following conditions, where $J_{-1}=P_{-1}=\emptyset$:

\medskip
(b$'$) the inclusion $P_{i-1}\cup f(J_{i-1})\,\i\,P_i$ extends to a map
$H_i\:M_i\to P_i$ of the mapping cylinder $M_i$ of the composition
$P_{i-1}\cup f(J_{i-1})\,\i\,P_i@>\phi_i>>Q_i$ so that
$(\phi_iH_i)^{-1}(C)=R_i^{-1}(C)$ for each $C\in K_i^*$, where
$R_i\:M_i\to Q_i$ is the projection.
\smallskip

(c$'$) the inclusion $J_{i-1}\cup f^{-1}(P_i)\,\i\,J_i$ extends to a map
$h_i\:N_i\to J_i$ of the mapping cylinder $N_i$ of the composition
$J_{i-1}\cup f^{-1}(P_i)\,\i\,J_i@>\psi_i>>\tilde Q_i$ so that
$(\psi_ih_i)^{-1}(C)=r_i^{-1}(C)$ for each $C\in\tilde K_i^*$, where
$r_i\:N_i\to\tilde Q_i$ is the projection.
\medskip

$X$ is said to {\it $k$-quasi-embed} in $Y$ if there exists a
$k$-quasi-embedding $f\:X\to Y$.
Embeddings $g,h\:X\emb Y$ (or more generally $k$-quasi-embeddings) are called
{\it $k$-quasi-concordant} if they can be joined by a $k$-quasi-embedding
$F\:X\x I\to Y\x I$.
Finally, a {\it $k$-quasi-isotopy} is a level preserving
$k$-quasi-concordance $F$ for which $pt$ can be replaced by $I$ in
the above diagram (that is, the composition $S_F\,\i\,X\x I\to I$ factors
through $P_0$ and all the $Q_i$'s).
\endremark

\remark{(iii)}
The above definitions are modelled after the controlled version of
the Penrose--Whitehead--Zeeman--Irwin trick as it appears in the Homma--Bryant
proof \cite{Br} of the Chernavskij--Miller Theorem (compare \cite{Me1; \S3})
that topological embeddings between PL manifolds in codimension $\ge 3$ are
approximable by PL ones.
In more detail, let us say that a (strong) $k$-quasi-embedding $f\:X\to Y$ is
{\it $\eps$-controlled} if $\phi_k^{-1}(C)$ and $f(\psi_k^{-1}(D))$ have
diameters $<\eps$ for each $C\in K_k^*$ and $D\in\tilde K_k^*$.

In these terms, Bryant's arguments show that, firstly, if $\xi\:X^n\to Y^m$
is a topological embedding between PL manifolds with $m-n\ge 3$, for each
$k\in\N$ there exists an $\eps=\eps(k,\xi)>0$ such that any generic PL
$\eps$-approximation $f$ of $\xi$ is a $c\eps$-controlled strong
$k$-quasi-embedding for some constant $c=c(m,n)$.
Secondly, the very same arguments of \cite{Br} prove that there exists
a constant $k(m,n)$ such that every $\eps$-controlled
$k(m,n)$-quasi-embedding $X\to Y$ is $\eps$-close to a PL embedding.

In particular, $k(m,n)$-quasi-embeddability of $X$ in $Y$ implies
PL embeddability, and similarly $k(m,n)$-quasi-concordance of PL embeddings
$X\emb Y$ implies PL concordance, which in turn implies (cf\. \cite{Me1})
ambient PL isotopy.

The higher-dimensional version of $k$-quasi-isotopy also seems to be of some
interest in connection with configuration spaces of $(k+2)$-tuples of
points (compare \cite{Me1; Conjecture 1.9}).
However, we will not pursue it any further in this paper.
\endremark

\remark{(iv)} There are some reasons to regard the Homma--Bryant technique as
somewhat complementary or alternative to the technique of gropes.
Indeed, the latter was first employed by Shtan'ko to prove a version of the
Chernavskij--Miller Theorem for compacta (cf\. \cite{Ed2}), and according to
\cite{Ed2; last line on p\. 96} it works to give yet another proof of the
Chernavskij--Miller Theorem (but note that, though Shtan'ko's Theorem implies
the Chernavskij--Miller Theorem via a theorem of Bryant--Seebeck, this does not
yield the desired proof, since the latter rests on the Chernavskij--Miller
Theorem itself).
Later, both gropes and Casson handles have played an important role in
$4$-manifold topology.
Now, both viewpoints  --- gropes and $k$-quasi-isotopy --- yield geometric
approaches to finite type invariants, respectively, of knots \cite{CT} and of
links ``modulo knots'' (see a further paper by the first author).
\endremark

\remark{(v)} In looking for an equivalence relation on classical PL links
such that any two PL links sufficiently close to a given topological link
are equivalent, one is rather naturally led to $k$-quasi-isotopy.
Recall that in codimension $\ge 3$, it was proved by Edwards \cite{Ed1}
that any two PL embeddings of a compact polyhedron $X$ into a PL manifold $Y$
sufficiently close to a given topological embedding are joined by a small
PL ambient isotopy.
So we could obtain the desired equivalence relation by examining general
position properties, essential for Edwards' Theorem.
Edwards' proof rests mainly on radial engulfing and on his ``Slicing Lemma'',
proved by a version of the Penrose--Whitehead--Zeeman--Irwin trick where
the usual application of engulfing is excluded.
An alternative proof of Edwards' Theorem in the case where $X$ is a PL
manifold can be given along the lines of \cite{Ze2; proof of Theorem 24}
using the Homma--Bryant technique in place of the
Penrose--Whitehead--Zeeman--Irwin trick.
(This is a more direct argument than the straightforward combination of
the relative version of the Chernavskij--Miller Theorem and the ``small
concordance implies small ambient isotopy'' theorem of \cite{Me1}.)
In either case, we arrive at some variant of $k$-quasi-isotopy.
\endremark

\head 2. Invariants of $1$-quasi-isotopy \endhead

\subhead \oldnos2.\oldnos1. First examples \endsubhead
A counterpart to the $0$-quasi-isotopically nontrivial Hopf link $\H$
(Fig\. 3a), the (left handed) Whitehead link $\W$ (Fig\. 3b) is not weakly
$1$-quasi-isotopic to the unlink.
Indeed, the crossing change formula \cite{Jin} for the Sato--Levine invariant
$\beta=\bar\mu(1122)$ of $2$-component links with zero linking number implies
its invariance under weak $1$-quasi-isotopy:
$$\beta(L_+)-\beta(L_-)=-n^2,\tag{2.1}$$ where $L_+,L_-\:S^1\sqcup S^1\emb S^3$
differ by a single positive%
\footnote{We call the process of self-intersection {\it positive} if, given
a choice of orientation of the component that crosses itself and a choice of
ordering of the two singular points $p,q$, the triple consisting of
the positive tangent vector at $L_-(p)$, of that at $L_-(q)$, and of the vector
$L_-(p)-L_-(q)$, agrees with a fixed orientation of the ambient space
(this is clearly independent of the two choices).}
self-intersection of one component, leading to two lobes whose linking numbers
with the other component are $n$ and $-n$.
The formula (2.1) along with the requirement that $\beta($unlink$)=0$ clearly
suffices to evaluate $\beta$ on any null-homotopic $2$-component link, so
$\beta(\W)=1$ and the assertion follows.

\fig 3

Next, one could ask, whether the three versions of $1$-quasi-isotopy are
distinct.
The untwisted left handed Whitehead double of $\W$, denoted $\W_2$
(Fig\. 3c), is clearly $1$-quasi-isotopically trivial but {\sl seemingly} not
strongly $1$-quasi-isotopically trivial since the Whitehead curve in
the solid torus cannot be engulfed (i.e\. $\W$ is not a split link).
Whereas a link of two trefoils $\W_2'$ obtained by a twisted
Whitehead doubling of each component of the Hopf link (Fig\. 3d) can be seen to
weakly $1$-quasi-isotop onto the unlink, but {\sl conjecturally} is nontrivial
up to $1$-quasi-isotopy.

\proclaim{Problem 2.1} Find invariants distinguishing the three versions of
$1$-quasi-isotopy.
\endproclaim

\subhead \oldnos2.\oldnos2. The generalized Sato--Levine invariant \endsubhead
Apart from these matters, let us see the difference between (either version of)
$1$-quasi-isotopy and (either version of) concordance.
Consider the well-known link $\M$ depicted in Fig\. 4a.
We call it the {\it Mazur link} (compare \cite{Ma}) since it was employed,
{\it inter alia,} in the construction of Mazur's contractible $4$-manifold.
The Mazur link is obviously F-isotopic%
\footnote{Two links are called {\it F-isotopic} if one can be obtained from
another by a sequence of substitutions of a PL-embedded component $K_i$ with
an arbitrary knot $K_i'$ lying in its regular neighborhood $V(K_i)$ and
homotopic to $K_i$ in $V(K_i)$, and of reverse operations.}
to the Hopf link and hence topologically
I-equivalent to it, by virtue of Giffen's shift-spinning construction
(see \cite{Hi; Theorem 1.9}).
It is known that $\M$ is not PL I-equivalent to $\H$ \cite{Rol}, but
nevertheless the {\it fake Mazur link} (cf\. \cite{Ma; p\. 470}) $\M'$,
shown in Fig\. 4b, is even PL concordant (i.e\. locally flat PL I-equivalent)
to $\H$.
This is clear from Fig\. 4c, where the ambient isotopy class of $\M'$ is
represented in a different way.

\fig 4

On the other hand, both $\M$ and $\M'$ are not weakly $1$-quasi-isotopic
to $\H$, as detected by the {\it generalized Sato--Levine invariant}
$\tl\beta$ (compare \cite{KL1; p\. 1351} and \cite{NS; Fig\. 6}).
By definition, $\tl\beta(L)=a_3(L)-a_1(L)(a_2(K_1)+a_2(K_2))$, where
$K_1$ and $K_2$ are the components of the link $L$, and $a_i$ denotes
the coefficient at $z^i$ of the Conway link polynomial.
It is easy to verify (cf\. \cite{Liv}), using the skein relation for
the Conway polynomial, that, given two links $L_+,L_-\:S^1\sqcup S^1\emb S^3$,
related by a single positive crossing change on one component so that the two
lobes' linking numbers with the other component are $n$ and $l-n$, where
$l=\lk(L_+)=\lk(L_-)$, we have
$$\tl\beta(L_+)-\tl\beta(L_-)=n(l-n).\tag{$2.2$}$$

Hence $\tl\beta(\H)\neq\tl\beta(\M)=\tl\beta(\M')$, meanwhile (2.2)
immediately implies invariance of $\tl\beta$ under weak $1$-quasi-isotopy.
This proves the case of links in the following

\proclaim{Theorem 2.2} PL concordance does not imply weak $1$-quasi-isotopy,
for both links and string links.
\endproclaim

Proceeding to the case of string links, let us first introduce some notation
to be used in the rest of this section.

If $L\:I_+\sqcup I_-\emb\R^2\x I_0$ is a string link, let $\rho L$ denote
the {\it reflection} $R\circ L\circ r$ of $L$, where $r$ and $R$ are the
self-homeomorphisms of $I_+\sqcup I_-$ and $\R^2\x I_0$ given by
the reflections of $I_+$, $I_-$ and $I_0$ in their midpoints.

Next, let $H_n$ denote the string link that can be described as the braid
$(\sigma_{12})^{2n}$ or the integer tangle $2n$ (see \cite{GK}).
The link $H:=H_1$ will be called the {\it Hopf string link} since its
numerator closure (see \S1) $numer(H)$ is the Hopf link $\H$.
The {\it linking number $\l$ closure} of a $2$-component string link $L$ is
defined to be the link $numer_l(L):=numer(L\#H_{l-n})$, where
$n=\lk(numer(L))$.

Finally, recall from \S1 that $W$ denotes the Whitehead string link, shown
three times in Fig\. 1; it can also be described as the rational tangle
$\frac1{1+\frac12}$.

\demo{Proof of the case of string links in Theorem 2.2}
The link $numer_1(W\#\rho W)$ is the fake Mazur link (Fig\. 4c), which is not
$1$-quasi-isotopic to the Hopf link.
Therefore the string link $W$ is not cancelled by $\rho W$ up to
$1$-quasi-isotopy.

However, $W\#\rho W$ is concordant to the trivial string link.
Indeed, for any string link $L\:I_+\sqcup I_-\emb\R^2\x I_0$, the embedding
$L\x\bold 1_I\:(I_+\sqcup I_-)\x I\to\R^2\x I_0\x I$ is a collection of
slicing disks for the restriction of $L\x\bold 1_I$ to
$(I_+\sqcup I_-)\x\partial I\cup (\{1\}\sqcup\{1\})\x I$, which is a string
link ambient isotopic to $L\#\rho L$.
\qed
\enddemo

\remark{Remarks. (i)}
An axiomatic definition of $\tl\beta$ can be given by the formula (2.2)
along with the requirement that $\tl\beta(\H_n)=0$ for each $n\in\Z$,
where $\H_n=numer(H_n)$.
This definition was used in \cite{KL1}, where $\tl\beta$ was related to the
Casson--Walker invariant of rational homology $3$-spheres, in \cite{AR},
where a direct proof of its existence was given, without references to
skein theory, and in \cite{Liv}.
An invariant defined by Gauss diagrams in \cite{PV} (the misprint in
the definition was corrected in \cite{\"Os} and \cite{AMR}) turned out to
coincide with $\tl\beta$ \cite{AMR} (as defined above and in \cite{KL1};
the axiomatic definition of $\tl\beta$ in \cite{AMR} uses a slightly different
choice of $\H_n$, discussed in \cite{KL1; p\. 1337}).
Yet another choice of $\H_n$ has to be made for the version of $\tilde\beta$
in \cite{Tr; \S10}.
Nakanishi and Ohyama recently proved that $2$-component links have the same
$\lk$ and $\tl\beta$ iff they are related by a sequence of the so-called
self-$\Delta$-moves \cite{NO}.
Relevance of $\tl\beta$ for magnetohydrodynamics is discussed in \cite{Ah}.
\endremark

\remark{(ii)} It follows from \cite{Co3; \S4} (see also \cite{Co2; \S9}) that
$\tl\beta$ is an integral lifting of Milnor's invariant $\bar\mu(1122)$.
To add some clarity to the subject (cf\. \cite{KL1; Addenda}), we give
a direct proof of this fact in Appendix.
\endremark
\medskip

By the Stallings Theorem, the quotients of the fundamental groups of
topologically I-equivalent links (or string links) by any term of the lower
central series are isomorphic \cite{St1}, \cite{Ca1}, \cite{Co1}, and
the isomorphism preserves the image of the peripheral structure in
the quotient (see \cite{Ca1; Remark 1}).
Thus we have

\proclaim{Corollary 2.3} No nilpotent quotient of the fundamental group, even
if equipped with the peripheral structure, can classify two-component links,
or string links, up to $1$-quasi-isotopy.
\endproclaim

\subhead \oldnos2.\oldnos3. The $\eta$-function and the Jin suspension
\endsubhead
Let us recall the definition of the $\eta$-function of Kojima and Yamasaki
\cite{KY}.
Given a linking number zero link $L(2S^1)=K_+\sqcup K_-\i S^3$, let $\tl K_+$
be a lift of $K_+$ into the infinite cyclic cover $X_-$ of $S^3\but K_-$, and
let $\tl K^*_+$ be the nearby lift of the zero pushoff of $K_+$ (which may be
a nonzero pushoff of $\tl K_+$ if $K_-$ happens to be unknotted).
Let $f(t)$ be any Laurent polynomial annihilating the class of $\tl K_+$ in
the $\Z[t^{\pm 1}]$-module $H_1(X_-)$,%
\footnote{The Wang sequence for a map $S^3\but K_-\to S^1$ representing
a generator of $H^1(S^3\but K_-)$, which is the same thing as the long
sequence associated to the sequence $0\to\Z[\Z]@>t-1>>\Z[\Z]@>t=1>>\Z\to 0$
of coefficient modules for $H_*(S^3\but K_-;\cdot)$, shows that multiplication
by $t-1$ is an isomorphism on $H_1(X_-)$.
Since $H_1(X_-)$ is finitely generated over $\Z[\Z]$, which is Noetherian,
$H_1(X_-)$ is Noetherian, therefore infinite divisibility by $t-1$ yields,
for each $m\in H_1(X_-)$, a Laurent polynomial $p$ such that $(t-1)^{-1}m=pm$.
It follows that $H_1(X_-)$ is torsion over $\Z[\Z]$.}
for example the Alexander polynomial of $K_-$.
Then the $1$-cycle $f(t)\tl K_+$ bounds a chain $\zeta$ in $X_-$.
Define
$$\eta^+_L(t)=\frac1{f(t)}\sum_{n=-\infty}^\infty
(\zeta\cdot t^n\tl K_+^*)t^n,$$
where $(\,\cdot\,)$ stands for the intersection pairing in $X_-$, and only
finitely many of the summands may be nonzero; thus $\eta^+_L(t)\in\Q(t)$,
the field of fractions of $\Z[t^{\pm 1}]$.
Then $\eta^+_L$ does not depend on the choice of $\zeta$ since $H_2(X_-)=0$,%
\footnote{The same argument as in the preceding footnote shows that
$H_2(X_-)$ is $\Z[\Z]$-torsion.
But $S^3\but K_-$ is homotopically $2$-dimensional, so
$H_2(X_-)=\ker\left[C_2(X_-)\to C_1(X_-)\right]$ is free over $\Z[\Z]$.}
and consequently on the choice of $f(t)$, for $\eta^+_L$ clearly remains
unchanged if $f(t)$ and $\zeta$ are simultaneously multiplied by some
polynomial $g(t)$.
By interchanging $K_+$ and $K_-$, one can similarly define $\eta^-_L$.
Since the projection of $\zeta$ to $S^3\but K_-$ has zero intersection number
with $K_+^*$, we have $\eta^+_L(1)=0$; therefore
$\eta^+_L(t)=(1-\eps)\left[f(t)^{-1}\sum_{n\ne 0}
(\zeta\cdot t^n\tl K_+^*)t^n\right]$,
where $\eps\:\Z[\Z]\to\Z\i\Z[\Z]$ is the augmentation, $\eps(t)=1$.

\remark{Remark}
If $F$ is a Seifert surface for $K_-$, disjoint from $K_+$, one can similarly
define a pairing $\left<\cdot\right>$ on disjoint cycles in a lift of
$S^3\but F$ into $X_-$ such that $\eta^+_L=\left<\tl K_+\cdot\tl K^*_+\right>$
and which descends to the Blanchfield pairing (see \cite{Hi}) on $H_1(X_-)$,
well-defined up to addition of (finite) polynomials.
The pairing $\left<\cdot\right>$ was used in \cite{Co2; \S7} to show that
the rational function $\eta^+_L(t)$ can be expanded as a rational power
series $\sum_{k=1}^\infty\beta^k_+z^k$ in the variable $z=-(t+t^{-1}-2)$,
where the coefficients $\beta^k_+$ are known as Cochran's derived
invariants (in particular, $\beta^1_+$ is the Sato--Levine invariant
$\beta$).
\endremark
\medskip

If $K_-$ happens to be unknotted, $\eta^+_L$ is by definition the Laurent
polynomial $\sum_{n\ge 1}\lk(\tl K_+,t^n\tl K_+)(t^n+t^{-n}-2)$, where
$\lk$ is the usual linking number in $\R^3$.
More generally, by a well-known unpublished observation of J. T. Jin,
the $\eta$-function satisfies the following crossing change formula for
a positive self-intersection of $K_+$:
$$\eta^+_{L'}-\eta^+_L=t^s+t^{-s}-2,\tag{2.3}$$
where $s$ is the linking number between $K_-$ and one of the lobes of
the singular component.
Indeed, we may assume that the trace of the homotopy between $K_+$ and $K_+'$
is a small disk $D$ such that $\partial D=K_+'-K_+$ as chains, and which meets
the common arc of $K_+$ and $K_+'$ transversally in an interior point $p$.
Then $f(t)\tl K'_+$ bounds the chain $\zeta'=\zeta+f(t)\tl D$, where
$\tl D$ is the appropriate lift of $D$, and we may assume that $\zeta'$ meets
the union of all translates of the lift $\tl D^*$ of the pushoff of $\tl D$
along the $1$-chain $f(t)\tl J^*$ where $\tl J^*$ is the lift of the pushoff
of an arc $J\i D$ joining $p$ with a point of $\partial D$.
We have
$$\eta^+_{L'}(t)-\eta^+_L(t)=\sum_{n\ne 0}(\tl D\cdot t^n\tl K_+^*)(t^n-1)+
(1-\eps)\left[\frac1{f(t)}\sum_{n\ne 0}(\zeta'\cdot t^n\partial\tl D^*)t^n
\right],$$
where the first summand of the right hand side equals $t^s+t^{-s}-2$ and
the second $(1-\eps)[1]=0$.
See \cite{JL..} for a different proof and \cite{KL2} for generalizations of
the formula (2.3).

Since the right hand side of (2.3) turns to zero when $s=0$, we get

\proclaim{Proposition 2.4} $\eta^+_L$ is invariant under weak
$1$-quasi-isotopy with support in $K_+$.
\endproclaim

A link $S^1\sqcup S^1\emb S^3$ is called {\it semi-contractible} \cite{Ki},
\cite{Ko} if each component is null-homotopic in the complement to the other.
In the class of semi-contractible links there is yet another way to think of
Cochran's invariants.
The {\it Jin suspension} \cite{Ki} of a semi-contractible link $L$ is
a link map $S^2\sqcup S^2\to S^4$, formed by the tracks of some
null-homotopies of the components of $L$ in each other's complement.
By the Sphere Theorem of Papakyriakopoulos, the result is well-defined up to
link homotopy (cf\. \cite{Ko}).
This construction, applied by Fenn and Rolfsen to the Whitehead link $\Cal W$,
produced the first example of a link map $S^2\sqcup S^2\to S^4$, nontrivial
up to link homotopy.

To define Kirk's invariant of a link map $\Lambda\:S^2_+\sqcup S^2_-\to S^4$,
we assume by general position that $f$ only has transversal double points.
Let $z=\Lambda(x)=\Lambda(y)$ be one, with $x,y\in S^2_+$.
Let $p$ be a path in $S^2_+$ joining $x$ to $y$, then $\Lambda(p)$ is
a closed loop in $S^4\but\Lambda(S^2_-)$.
Denote by $n_z$ the absolute value of the linking number between this loop
and $\Lambda(S^2_-)$, and by $\eps_z$ the sign, relating the orientations of
the two sheets crossing at $z$ to that of $S^4$.
Then $\sigma_+(t):=\sum_z \eps_z(t^{n_z}-1)$, together with the analogously
defined $\sigma_-(t)$, form Kirk's invariant $\sigma(\Lambda)$ of link
homotopy of $\Lambda$ \cite{Ki}.
By (2.3), the pair $(\eta^+_L,\eta^-_L)$ for a semi-contractible link $L$
contains the same information as the Kirk invariant of the Jin suspension
of $L$.
So Proposition 2.4 implies

\proclaim{Proposition 2.5}
Let $L$ and $L'$ be semi-contractible links that can be joined by two weak
$1$-quasi-isotopies such that the $i^{\text{th}}$ carries the $i^{\text{th}}$
component by a locally-flat PL isotopy.
Then the Jin suspensions of $L$ and $L'$ have identical Kirk's invariant.
\endproclaim

\remark{Remark}
According to P. Teichner, it remains an open problem (as of Spring 2003)
whether two link maps $S^2\sqcup S^2\to S^4$ with the same Kirk invariant
may be not link homotopic (compare \cite{BT}).
\endremark

\subhead \oldnos2.\oldnos4. Koschorke's refined Kirk invariant \endsubhead
We conclude this section with a result on link maps, imposing a constraint
on existence of invariants of $1$-quasi-isotopy of string links.
(The constraint is discussed in Remark (i) at the end of the section.)

Under a {\it fibered disk link map}
$\Lambda\:(I_+\sqcup I_-)\x I\to(\R^2\x I_0)\x I$ we understand a
self-link-homotopy of the trivial string link.
(In terms of the map $\Lambda$, we assume that
$\Lambda(I_\pm\x\{t\})\i\R^2\x I_0\x\{t\}$ for each $t\in I$, and that
$\Lambda(\pm,s,t)=(\pm p,s,t)$ whenever $(s,t)\in\partial (I\x I)$,
where $p\in\R^2\but\{0\}$ is a fixed point.)
The {\it trivial} fibered disk link map is given by
$\Lambda_0(\pm,s,t)=(\pm p,s,t)$ for all $s$ and $t$.
To every fibered disk link map $\Lambda$ one can associate a spherical
link map $Numer(\Lambda)\:S^2_+\sqcup S^2_-\to S^4$, defined to be the
composition of $\Lambda/\partial\:(I_+\x I)/\partial\sqcup
(I_-\x I)/\partial\to\R^2\x ((I_0\x I)/\partial)$ and the inclusion
$\R^2\x S^2\i D^2\x S^2\cup S^1\x D^3=\partial (D^2\x D^3)$.

For fibered disk link maps there is a refined version of the Kirk
invariant, well-defined up to {\sl fibered} link homotopy \cite{Ko}.
Given a $\Lambda$ as above, by general position we may assume that it only
has transversal double points.
Let $z=\Lambda(x,t)=\Lambda(y,t)$ be one, with $x,y\in I_+$.
There is a path $p$ in $I_+$ joining $x$ to $y$, and it is unique up to
homotopy relative to the endpoints.
Then $\Lambda(p,t)$ is a closed loop in $\R^2\x I_0\x I\but\Lambda(I_-\x I)$,
which has a well-defined linking number with $\Lambda(I_-\x I)$.
More precisely, let $l_z$ denote the algebraic number of intersections
between $\Lambda(I_-\x\{t\})$ and a generic disk, spanned by $\Lambda(p,t)$ in
$\R^2\x I_0\x\{t\}$.
We emphasize that not only the absolute value, but also the sign of $l_z$
is well-defined, due to the natural ordering of $x$ and $y$ in $I_+$.
Denote by $\eps_z$ the sign, relating the orientations of the two sheets
crossing at $z$ to that of the ambient $4$-space.
We define $\tl\sigma_+$ to be the Laurent polynomial
$\sum_z\eps_z(t^{l_z}-1)$, where $z$ runs over the double points of
the $+$-component $\Lambda|_{I_+\x I}$, and $\tl\sigma_-$ by interchanging
$I_+$ and $I_-$.

To see that $\tl\sigma:=(\tl\sigma_+,\tl\sigma_-)$ is unchanged under
a fibered link homotopy between fibered disk link maps $\Lambda$ and
$\Lambda'$, approximate $\Lambda'$ arbitrarily closely by a $\Lambda''$,
fibered {\sl regularly} homotopic with $\Lambda$ (cf\. \cite{RS; 4.4}).
Clearly, $\tl\sigma(\Lambda'')=\tl\sigma(\Lambda')$, and since $l_z$ is
a continuous function on the $1$-manifold of double points of the fibered
regular homotopy, $\tl\sigma(\Lambda'')=\tl\sigma(\Lambda)$.
Of course, $\sigma(Numer(\Lambda))=(\phi\oplus\phi)\tl\sigma(\Lambda)$, where
$\phi\:\Z[t^{\pm 1}]\to\Z[t]$ denotes the $\Z$-homomorphism given by
$t^n\mapsto t^{|n|}$.

Here are some examples.
Semi-contractible string links and their Jin suspensions, denoted
$\frak J(L)$, are defined in the obvious fashion.
Let $h^\pm_n\:S^1_+\sqcup S^1_-\to S^1_+\sqcup S^1_-$ be a map that has
degree $n$ on $S^1_\pm$ and degree $1$ on $S^1_\mp$.
Let us fix some $2$-component string link $W_\pm^n$, whose linking number zero
closure $numer_0(W_\pm^n)$ is isotopic to a $C^0$-approximation of
the composition of $h^\pm_n$ and the Whitehead link $\Cal W$.
In particular, $W_+^1$ and $W_-^1$ may be chosen to coincide with
the Whitehead string link $W$ (shown three times in Fig\. 1).
It is not hard to see (cf\. \cite{Ko}) that
$$\tl\sigma(\frak J(W_+^n))=(\tfrac{n^2\!\!\!}2\ (t+t^{-1}-2)+
\tfrac n2(t-t^{-1}),1-t^n),$$ and it follows that
$$\tl\sigma(\frak J(W_+^n\#\rho W_+^n))=(n(t-t^{-1}),t^{-n}-t^n).$$

\proclaim{Theorem 2.6}
(a) \cite{Ki} $\im\sigma=\ker\delta$, where
$\delta\:\Z[t]\oplus\Z[t]\to\Z\oplus\Z\oplus\Z$ is given by
$$(f,g)\mapsto (f|_{t=1},\ g|_{t=1},\ f'+g'+f''+g''|_{t=1}).$$

(b) $\im\tl\sigma=\ker\tl\delta$, where
$\tl\delta\:\Z[t^{\pm 1}]\oplus\Z[t^{\pm 1}]\to\Z\oplus\Z\oplus\Z\oplus\Z$
is given by
$$(f,g)\mapsto (f|_{t=1},\ g|_{t=1},\ f'+g'|_{t=1},\ f''+g''|_{t=1}).$$
\endproclaim

\demo{Proof of (b)}
By the definition of $\tl\sigma$, it satisfies
$\tl\sigma_+|_{t=1}=\tl\sigma_-|_{t=1}=0$.
Let $h_t$ be a self-link-homotopy of the trivial string link.
Then $numer_0(h_t)$ is a self-link-homotopy of the unlink, and $numer_1(h_t)$
is a self-link-homotopy of the Hopf link.
Since the generalized Sato--Levine invariant is a well-defined invariant
of links, it does not change under each homotopy.
On the other hand, it jumps by
$\sum_z\eps_zl_z(0-l_z)=-\tl\sigma_0'-\tl\sigma_0''|_{t=1}$ in
the former case, and by $\sum_z\eps_z l_z(1-l_z)=-\tl\sigma_0''|_{t=1}$ in
the latter, where $\tl\sigma_0=\tl\sigma_++\tl\sigma_-$.
Thus $\im\tl\sigma\i\ker\tl\delta$.

The reverse inclusion is due to Koschorke \cite{Ko}, and follows easily
from the above mentioned facts.
Indeed,
$\sigma(Numer(\frak J(W_\pm^n)))=(\phi\oplus\phi)\tl\sigma(\frak J(W_\pm^n))$
generate $\ker\delta$.
On the other hand, $\tl\sigma(\frak J(W_\pm^n\#\rho W_\pm^n))$
generate the kernel of the $\Z$-homomorphism
$\Z[t^{\pm 1}]\oplus\Z[t^{\pm 1}]\to\Z[t]\oplus\Z[t]\oplus\Z$, given by
$(f,g)\mapsto (\phi(f),\phi(g),f'+g'|_{t=1})$. \qed
\enddemo

The first part of the above argument yields

\proclaim{Proposition 2.7} If two semi-contractible $2$-component links
(string links) are weakly $1$-quasi-isotopic, their Jin suspensions have
identical $\sigma_\pm'(1)+\sigma_\pm''(1)$
(resp\. $\tl\sigma_\pm'(1)$ and $\tl\sigma_\pm''(1)$).
\endproclaim

The following seems to be implicit in \cite{Ko}.

\proclaim{Theorem 2.8} Fibered concordance does not imply fibered link
homotopy.
\endproclaim

\demo{Proof} Consider the fibered disk link map $\frak J(W\#\rho W)$.
By a standard construction (cf\. the proof of Theorem 2.2) it is fiberwise
concordant with the trivial fibered disk link map.
However, they are not fiberwise link homotopic, since
$\tl\sigma(\frak J(W\#\rho W))=(t-t^{-1},t^{-1}-t)\ne(0,0)$. \qed
\enddemo

This yields a counterexample (for disk link maps) to the 1-parametric version
of Lin's Theorem \cite{Lin} (see also \cite{Gi},\cite{Go}) that singular
concordance implies link homotopy in the classical dimension, and to
the fibered version of Teichner's Theorem \cite{BT; Theorem ~5}.

\remark{Remarks. (i)} The definition of $\tl\sigma$ applies equally well to
an arbitrary link homotopy of string links, rather than a self-link-homotopy
of the trivial string link; we denote the extended invariant by
$\bar{\tl\sigma}=(\bar{\tl\sigma}_+,\bar{\tl\sigma}_-)$.
Following the proof of ``$\i$'' in Theorem 2.6(b), the invariants
$\beta(numer_0(L))$ and $\tl\beta(numer_1(L))$ of a string link $L$
can now be viewed as determined by the extended invariant
$\bar{\tl\sigma}_0:=\bar{\tl\sigma}_++\bar{\tl\sigma}_-$ of any link homotopy
$H$ from $L$ to the trivial string link, namely
$$\beta(numer_0(L))=-\bar{\tl\sigma}_0'(H)-\bar{\tl\sigma}_0''(H)|_{t=1}
\qquad\text{and}\qquad\tl\beta(numer_1(L))=-\bar{\tl\sigma}_0''(H)|_{t=1}.$$
On the other hand, the ``$\supset$'' part of Theorem 2.6(b) implies that
any invariant of a string link $L$ that is a $\Z$-linear function of
$\bar{\tl\sigma}$ of a link homotopy from $L$ to the trivial string link, is
a linear combination of $\beta\circ numer_0$ and $\tl\beta\circ numer_1$.
\endremark

\remark{(ii)} The $\bmod 2$ reduction of $\bar{\tl\sigma}_+'|_{t=1}$, i.e\.
$\sum_z\eps_z l_z\pmod 2$, where $z$ runs over the double points of
the $+$-component, is nothing but Hudson's obstruction \cite{Hud} to
embedding the $+$-component by homotopy in the complement of the
$-$-component.
Note that the $\bmod 2$ reduction of $\bar{\tl\sigma}_+''|_{t=1}$ is
identically zero.
\endremark

\remark{(iii)} The string links $W^{-1}_+$ and $W^{-1}_-$ may be chosen to
coincide with the string link $W^{-1}$ that can be described as the rational
tangle $\frac1{2+\frac12}$.
Its reflection $\bar W:=\rho W^{-1}$ may be called the right handed Whitehead
string link, since its linking number zero closure $numer_0(\bar W)$ is
the right handed Whitehead link $\bar\W$ (i.e\. $\W$ composed with a
reflection of $S^3$ and precomposed with a reflection of one component).
Interestingly, $W\#\bar W$ and $\bar W\# W$ have zero invariants
$\beta\circ numer_0$ and $\tl\beta\circ numer_1$ as well as
$\eta^+\circ numer_0$ and $\eta^-\circ numer_0$, but seem unlikely to be
weakly null-$1$-quasi-isotopic.
\endremark

\head 3. The fundamental group \endhead

\subhead \oldnos3.\oldnos1. Preliminaries \endsubhead
Let us fix some notation.
Unless otherwise mentioned, we work in the PL category.
For any link $L\:mS^1\emb S^3$ we write $\pi(L)$ for its fundamental group
$\pi_1(S^3\but L(mS^1))$.
If $G$ is a group, the conjugate $h^{-1}gh$ of $g\in G$ by $h\in G$ will be
denoted by $h^g$, and the normal closure of a subgroup $H$ in $G$ by $H^G$.
We also write $g^{\pm h}$ for $(g^h)^{\pm 1}$.

For the following lemma, suppose that links $L_+,L_-\:mS^1\emb S^3$
are related by a crossing change on the $i^{\text{th}}$ component, and
$L_s\:mS^1\to S^3$ is the intermediate singular link, with lobes $\ell$
and $\ell'$.
We may assume that the images of $L_+$, $L_-$ and $L_s$ are contained in
some $\R^3\i S^3$ and coincide outside the cube $Q:=[-1,1]^3$, meeting
$L_s(mS^1)$ in $[-1,1]\x\{0\}\x\{0\}\cup\{0\}\x\{0\}\x[-1,1]$.
Attach $1$-handles $H,H'\cong D^2\x I$ to the interiors of four faces
$F_\pm$, $F'_\pm$ of $Q$ so that $Q\cup H$, $Q\cup H'$ and
$N:=Q\cup H\cup H'$ are regular neighborhoods of $\ell$, $\ell'$ and
$L_s(S^1_i)$ in $S^3\but L_s(mS^1\but S^1_i)$.

By a {\it meridian} of the lobe $\ell$ we mean the image in
$\pi(L_s):=\pi_1(S^3\but N)$ of a generator $\alpha$ of
$\pi_1(H\cap\partial N)\simeq\Z$ under the inclusion induced
homomorphism, which is determined by a path $\rho$ joining the basepoints
of these spaces.
The {\it longitude} of $\ell$, corresponding to $\rho$, is the image in
$\pi(L_s)$ of a generator $\beta$ of
$\pi_1((H\cup F_+\cup F_-)\cap\partial N)$ that can be represented by a loop
crossing the edge $F_+\cap F_-$ of $Q$ (geometrically) exactly once, and
having zero linking number with $\ell$.

\proclaim{Lemma 3.1} Suppose that links $L_+$ and $L_-$ differ by
a single crossing change on the $i^{\text{th}}$ component so that
the intermediate singular link $L_s$ is a $k$-quasi-embedding.
Let $\ell=J_0$ be the lobe of $L_s$ as in the definition of
$k$-quasi-embedding.
Let $\mu,\lambda\in\pi(L_s)$ be a meridian of $\ell$ and
the corresponding longitude, and $\mu_+,\lambda_+$ be their images in
$\pi(L_+)$.
Then $$\lambda_+\in\left<\mu_+\right>_k^{\pi(L_+)},$$ where for a subgroup
$H$ of a group $G$ we set $H_0^G=G$ and inductively $H_{k+1}^G=H^{H_k^G}$.
\endproclaim

It may be easier to guess the following argument than to read it.

\demo{Proof} The proof of the general case is best exposed by analogy with
the cases $k=1,2$, considered below.

Without loss of generality, $Q\cup H\i P_1$, where $P_1$ is as
in the definition of $k$-quasi-embedding.
Let $L_s(p)=L_s(q)$ be the double point of $L_s$, and $y$ the basepoint
of $H\cap\partial N$.
Let $\rho$ be a path joining the basepoint of $S^3\but N$ to $y$, and
$o\:S^1\to H\cap\partial N$, $l\:S^1\to (H\cup F_+\cup F_-)\cap\partial N$ be
some representatives of $\alpha$, $\beta$ above.
Pick a generic PL null-homotopy $F\:D^2\to P_1$ of the loop $l$, transversal
to $L_s(S^1_i)$ and such that $F(\Int D^2)$ is disjoint from $Q$.
Suppose, for definiteness, that $l$ is clockwise with respect to a fixed
orientation of $D^2$.
Note that $y=F(y)$.

If $N$ is small enough, each intersection point $L_s(p_j)=F(q_j)$ corresponds
to a small circle $C_j\i F^{-1}(N)$, circling around $q_j$.
Pick any $y_j\in C_j$, and let $I_j\i D^2$ be an arc joining $y$ to $y_j$ and
disjoint from all $C_u$'s, except at $y_j$.
Let us parameterize $F(C_j)$ by a clockwise loop $l_j$ with ends at $F(y_j)$,
and parameterize $F(I_j)$ by a path $r_j$ directed from $F(y)$ to $F(y_j)$.
Then the loop $l$ is homotopic in $F(D^2)\but\Int N$ to
$\bar r_1l_1r_1\dots\bar r_nl_nr_n$, so $\lambda_+$ is a product of
the meridians $[\bar\rho\bar r_jl_jr_j\rho]$ to the $i^{\text{th}}$ component.
But any two meridians to the same component are conjugate, thus
$\lambda_+\in\left<\mu_+\right>^{\pi(L_+)}$.
This completes the proof in the case $k=1$.

For $k\ge 2$ we need an additional ingredient for our construction.
Pick a parallel push-off $h$ of $L_+(S^1_i)$ into $\partial N$ such that
$h(p)=y$.
Each $y_j$ could have been chosen so that $F(y_j)=h(p_j)$.

Now let $I'_j$ denote the subarc of $J_1$ joining $p$ to $p_j$, and let
$\tl r_j$ be a path parameterizing $h(I'_j)$ in the direction from
$h(p)=F(y)$ to $h(p_j)=F(y_j)$.
Pick a generic PL null-homotopy $F_j\:D^2_j\to P_2$ of the loop
$\bar r_j\tl r_j$, transversal to $L_s(S^1_i)$ and such that
$F_j(\Int D^2_j)$ is disjoint from $Q$.
For convenience, we assume that $\bar r_j\tl r_j$ is counter-clockwise with
respect to a fixed orientation of $D^2_j$.
Note that $F_j(y)=y$ and $F_j(y_j)=y_j$.

If $N$ is small enough, each intersection point
$L_s(p_{ju})=F_j(q_{ju})$ corresponds
to a small circle $C_{ju}\i F_j^{-1}(N)$, circling around $q_{ju}$.
Pick any $y_{ju}\in C_{ju}$, and let $I_{ju}\i D^2_j$ be an arc joining $y$ to
$y_{ju}$ and disjoint from all $C_{ju}$'s, except at $y_{ju}$.
Let us parameterize $F_j(C_{ju})$ by a clockwise loop $l_{ju}$ with ends at
$F_j(y_{ju})$, and parameterize $F_j(I_{ju})$ by a path $r_{ju}$ directed
from $F_j(y)$ to $F_j(y_{ju})$.

Then the loop $\bar r_j\tl r_j\tl l_j$ is null-homotopic in
$F_j(D^2_j)\but\Int N$, where $\tl l_j$ denotes
$\bar r_{j1}l_{j1}r_{j1}\dots\bar r_{jn_j}l_{jn_j}r_{jn_j}$.
It follows that the loop $\bar r_jl_jr_j$ is homotopic in
$F_j(D^2_j)\but\Int N$ to
$\Bar{\Tilde l}_j\Bar{\Tilde r}_jl_j\tl r_j\tl l_j$.
But the loop $\Bar{\Tilde r}_jl_j\tl r_j$ is homotopic in $N\but L_+(S^1_i)$
(though possibly not in $N\but L_s(S^1_i)$) to $o$, so
$\bar\rho\bar r_jl_jr_j\rho$ is homotopic in $S^3\but L_+(mS^1)$ to
$(\bar\rho\Bar{\Tilde l}_j\rho)\bar\rho o\rho(\bar\rho\tl l_j\rho)$, which
represents $\mu_+^{[\tl l_j]}$.
Thus $\lambda_+$ is a product of conjugates of $\mu_+$ by products of
the meridians
$[\bar\rho\bar r_{ju}l_{ju}r_{ju}\rho]\in\left<\mu_+\right>^{\pi(L_+)}$. \qed
\enddemo

\remark{Remark} To continue the proof for $k\ge 3$, we should have chosen each
$y_{ju}$ so that $F_j(y_{ju})=h(p_{ju})$.
\endremark
\medskip

\subhead \oldnos3.\oldnos2. Higher analogues of Milnor's group \endsubhead
Let $L$ be an $m$-component link and $M\i\pi(L)$ be the set of all meridians of
$L$, which is closed under inverse and conjugation, though is not a subgroup.
Let $\mu_k$ denote the subgroup of $\pi(L)$ generated by the set
$$\{[m,m^g]\mid m\in M, g\in\left<m\right>_k^{\pi(L)}\},$$
where $[x,y]=x^{-1}y^{-1}xy$.
This subgroup is normal since $[m,m^g]^h=[m^h,(m^h)^{g^h}]$.
Moreover, since any two meridians to one component are either conjugate or
inverse-conjugate, while the commutator identity $[x,y^{-1}]=[y,x]^{y^{-1}}$
implies the equality $[m^{-1},(m^{-1})^g]=[m,m^g]^{(m^gm)^{-1}}$, our $\mu_k$
coincides with the normal closure of
$$\{[m_i,m_i^g]\mid i=1,\dots,m;\ g\in\left<m_i\right>_k^{\pi(L)}\},$$
where $\{m_1,\dots,m_m\}$ is a fixed collection of meridians.
Recall that $\pi(L)/\mu_0$ is Milnor's group $\G(L)$, invariant under
homotopy \cite{Mi; Theorem ~2}, \cite{Le; Theorem ~4}.
In the same spirit we consider the quotient $\G_k(L)=\pi(L)/\mu_k$ (this
agrees with the notation of \S1, as we will see from Theorems 3.2 and 3.5).%
\footnote{It also follows from the proofs of Theorems 3.2 and 3.5 that the set
$\Lambda_i$ from the definition of the peripheral structure $\P_k$ in \S1
coincides with the right coset $N_k(m_i)\bar l_i$ of the subgroup $N_k(m_i)$
(defined in the statement of Theorem 3.2) containing the class $\bar l_i$ of
the longitude $l_i$ corresponding to a representative $m_i$ of $\bar m_i$.}

\proclaim{Theorem 3.2} If two links
$L,L'\:mS^1\emb S^3$ are $k$-quasi-isotopic then the groups
$\G_k(L)$, $\G_k(L')$ are isomorphic.

Moreover, the isomorphism takes any element $m$ representing a meridian of
$L(S^1_i)$ onto an element representing a meridian of $L'(S^1_i)$, whereas
the element representing the corresponding longitude of $L'(S^1_i)$ is
the image of the product of an element $\xi$ of
$$N_k(m)=\left<[g^{-1},g^m]\mid g\in\left<m\right>_k^{\G_k(L)}\right>$$ and
the element representing the corresponding longitude of $L(S^1_i)$.

Furthermore, $\G_k(L)$ is functorially invariant under $k$-quasi-isotopy of $L$
with respect to the quotient map $\pi(L)\to\pi(L)/\mu_k=\G_k(L)$.
\endproclaim

\demo{Proof}
It suffices to consider the case where $L'$ is obtained from $L$ by a
$k$-quasi-isotopy with a single self-intersection of the $i^{\text{th}}$
component.
We can assume that there is a plane projection, where this self-intersection
corresponds to replacement of an underpass by the overpass.
Say, an arc $L(I_-)$ underpasses an arc $L(I_+)$, and the arc $L'(I_-)$
overpasses the arc $L'(I_+)$, where $I_+$ and $I_-$ are disjoint arcs in
$S^1_i$.
Furthermore, we can assume that the $k$-quasi-isotopy has its support in
$I_+\cup I_-$.

Suppose that $\partial I_-=\{p_1,p_2\}$ and $\partial I_+=\{p_3,p_4\}$ so that
the points $p_1,p_2,p_3,p_4$ lie in the cyclic order on $S^1_i$.
For $j=1,2,3,4$ let $m_j$ be the meridian, starting and ending in the basepoint
(which we assume to lie above the plane) and underpassing $L'(S^1_i)$ once
below the point $L'(p_j)$ in the direction inherited from the main underpass.
The groups $\pi(L)$ and $\pi(L')$ can be presented by the same generators,
including (the classes of) $m_j$'s, and by the same relations, including
$m_1 m_2=m_4 m_3$, with one exception: the relation $m_1=m_3$ in $\pi(L)$ is
replaced with the relation $m_2=m_4$ in $\pi(L')$ (cf\. \cite{Mi},
\cite{Ca2; proof of Lemma 1 in Lecture 1}, \cite{Le; p\. 386}).

Let $L_s$ denote the singular link arising in the $k$-quasi-isotopy between
$L$ and $L'$, and let $\lambda$ denote the longitude of the lobe $J_0=\ell$
of $L_s(S^1_i)$.
Then $m_2=m_3^{\lambda}$ both in $\pi(L')$ and in $\pi(L)$.
By Lemma 3.1, $\lambda\in\left<m_2\right>_k^{\pi(L)}$ and analogously for
$\pi(L')$.
Therefore the commutator $[m_2,m_3]$ lies in $\mu_k\pi(L)$, and consequently
the relation $m_3 m_2=m_2 m_3$ holds in the quotient $\G_k(L)$.
Analogously, the same relation holds in $\G_k(L')$.
Now, this latter relation and either two of the relations $m_1=m_3$, $m_2=m_4$,
$m_1 m_2=m_4 m_3$ imply the third one, which means that all of these relations
do hold both in $\G_k(L)$ and in $\G_k(L')$.
Thus these groups are isomorphic and the first assertion is proved.
The assertion on the longitudes now follows from Lemma 3.3 below, since,
substituting $g=\lambda^{m^{-1}}$ we obtain $\xi=[g^m,g^{-1}]$. \qed
\enddemo

\proclaim{Lemma 3.3} \cite{Le; proof of Theorem 4 (note that there $[x,y]$
means $xyx^{-1}y^{-1}$)}
If $L$ and $L'$ are related by a homotopy with a single self-intersection of
the $i^{\text{th}}$ component, then, in the above notation,
$\xi=[\lambda,[\lambda^{-1},m^{-1}]]$.
\endproclaim

Let us call a singular link with one double point satisfying the conclusion
of Lemma 3.1 a {\it virtual $k$-quasi-embedding}.
While all our major results on $k$-quasi-isotopy clearly hold for virtual
$k$-quasi-isotopy (defined in the obvious way), coincidence of the two
notions is unclear for $k\ge 2$.

\proclaim{Problem 3.4} Does virtual $k$-quasi-isotopy imply $k$-quasi-isotopy,
for all $k$?
\endproclaim

\subhead \oldnos3.\oldnos3. Maximality and nilpotence of $\G_k(L)$ \endsubhead
In view of Theorem 3.2, the following implies, in particular, that no
difference between $k$-quasi-isotopy, virtual $k$-quasi-isotopy and strong
$k$-quasi-isotopy can be seen from functorially invariant (in the sense of \S1)
quotients of $\pi(L)$.

\proclaim{Theorem 3.5} If $\F(L)$ is a group, functorially invariant under
strong $k$-quasi-isotopy with respect to an epimorphism $q_L\:\pi(L)\to\F(L)$,
then $q_L$ factors through the quotient map $\pi(L)\to\pi(L)/\mu_k=\G_k(L)$
for every link $L$.
\endproclaim

\demo{Proof} In other words, we want to show that each relation imposed on
$\pi(L)$ in order to obtain $\G_k(L)$ must also be satisfied in any quotient
$\F(L)$, functorially invariant under strong $k$-quasi-isotopy.

We begin by noting that the requirement of functorial invariance can be
restated as follows.
In the notation of the definition of functorial invariance, let $p$ be a path
joining the double point of the singular link $L_s$ to the basepoint, and $R$
be a small regular neighborhood of $p$.
We may assume that each $L_\pm$ is sufficiently close to $L_s$, so that
$L_\pm(mS^1)$ meets $R$ in a couple of arcs properly embedded in $R$.
Now functorial invariance of $\F(L)$ is equivalent to the requirement that,
starting with an arbitrary allowable singular link $L_s$ and resolving it in
the two ways, we have for the obtained links $L_+$, $L_-$ that each of the
compositions $\pi_1(R\but L_\pm(mS^1))\overset i_*\to\to\pi(L_\pm)\to\F(L_\pm)$
is the abelianization $\Z*\Z\to\Z\oplus\Z$ (see \cite{Le; \S4}; compare
\cite{Mi; Fig\. ~2}).

Next, notice that there are two loops contained in $R$ and representing
meridians $m$, $m^g$ of $L_+$ conjugated by the element $g$ represented by the
loop $l$ starting at the basepoint, going along $\bar p$ (i.e\. the path $p$
in the reverse direction), then along either lobe $q$ of the singular component
of $L_s$ (which we assume disjoint from $L_+$ and $L_-$), and finally back
along $p$.
Conversely, given any link $L$, any path $p$ joining a point $x$ near a
component $K$ of $L$ to the basepoint, and any $g\in\pi(L)$, we can represent
$g$ as a loop $l=\bar pqp$ for some loop $q$ starting at $x$, and isotop $L$
by pushing a finger along $q$ to obtain a link $L_+$, and then push further
through a singular link $L_s$ to obtain a link $L_-$ as in the above situation.
\footnote{Notice that this construction depends on the choice of representative
$l$ of $g$; for instance, tying a local knot on the loop $q$ is likely to
change the isotopy class of $L_-$.}
Let $m=m(p)$ denote the meridian represented by the loop $\bar pop$ where $o$
is a small loop around $K$ starting at $x$.

Then it suffices to prove that for any link $L$ and any path $p$ as above, an
arbitrary element $g\in\left<m\right>_k^{\pi(L)}$, where $m=m(p)$, can be
represented by a loop $l=\bar pqp$ such that the corresponding singular link
$L_s$ (obtained using $L$ and $q$ as above) is allowable with respect to
$k$-quasi-isotopy.
The latter will follow once we show that $l\cup o\i B_1\i\dots\i B_k\i S^3\but
(L(mS^1)\but K)$ for some balls $B_1,\dots,B_k$ such that $L(mS^1)\cap B_i$ is
contained in an arc contained in $B_{i+1}$ (indeed, then without loss of
generality the whole finger will lie in $B_1$).

Since $g\in\left<m\right>_k^{\pi(L)}$, we can write $g$ as
$g=m^{\iota_1 g_1}\dots m^{\iota_n g_n}$ where each
$g_i\in\left<m\right>_{k-1}^{\pi(L)}$ and each $\iota_i=\pm 1$.
Analogously, each $g_i=m^{\iota_{i1}g_{i1}}\dots m^{\iota_{i,n_i}g_{i,n_i}}$
where $g_{ij}\in\left<m\right>_{k-2}^{\pi(L)}$ and $\iota_{ij}=\pm 1$, and
so on.
We will refer to these as $g_{\I}$'s and $\iota_{\I}$'s where $\I$
runs over a set $S_g$ of multi-indices of length $\le k$, including the
multi-index of length $0$.
Let $S_{g,i}$ be the subset of $S_g$ consisting of all multi-indices of length
$i$, let $\phi\:S_{g,i}\to S_{g,i-1}$ be the map forgetting the last index,
and let $n_{\I}=|\phi^{-1}(\I)|$.

The path $p$ possesses a small regular neighborhood
$T\i S^3\but(L(mS^1)\but K)$ containing the loop $o$ and such that $T\cap L$
is an arc.
Let us fix a PL homeomorphism $T\simeq [-1,1]\x I^2$, where $I=[0,1]$, such
that $T\cap L=\frac12\x\frac12\x I$, and $[-1,0]\x I^2$ is a small regular
neighborhood of the basepoint.
Let us represent the conjugators $g_\J$, $\J\in S_{g,k}$, of the last
stage by pairwise disjoint embedded paths $l_\J$ starting in
$[-\frac 1{k+1},0)\x I\x 1$, ending in $[-1,-\frac 1{k+1})\x I\x 1$ and
disjoint from $T$ elsewhere.
Also let us represent $m$ by $|S_g|$ distinct paths $r_{\I}$,
$\I\in S_g$, setting $r_{\I}=([0,\frac34]\x\{\frac14,\frac34\}\cup
\{\frac34\}\x[\frac14,\frac34])\x f(\I)$ with orientation given by
$\iota_{\I}$, where $f\:S_g\to [0,1]$ is an injective map such that
$f(S_{g,i})\i(\frac i{k+1},\frac{i+1}{k+1}]$ for each $i$.
Notice that the arc $r_{\I}$ is properly contained in the ball
$D_{\I}=I^2\x U$ where $U$ is a closed neighborhood of $f(\I)$ in
$I\but f(S_g\but\{\I\})$.

Assume that for some $i\in\{1,\dots,k\}$, each of the elements $g_{\I}$,
$\I\in S_{g,i}$ is represented by an embedded path $l_{\I}$ with the
initial endpoint in $[-1+\frac i{k+1},0)\x I\x \{\frac{i+1}{k+1}\}$ and the
terminal endpoint in $[-1,-1+\frac i{k+1})\x I\x\{\frac{i+1}{k+1}\}$, disjoint
from $[-1,1]\x I\x [0,\frac{i+1}{k+1}]$ except in the endpoints, and disjoint
from the other paths $l_\J$, $\J\in S_{g,i}$.
Let us use this to construct for each $\I\in S_{g,i-1}$ an embedded path
representing $g_{\I}$.
We have $g_{\I}=m^{\iota_1 g_{\I,1}}\dots m^{\iota_{n_{\I}}
g_{\I,n_{\I}}}$, where each $\iota_i=\pm 1$, and each of the
conjugators $g_\J$, $J\in\phi^{-1}(I)$, is represented by an embedded
path $l_\J$ by the assumption.
Let us connect the initial endpoint of each $l_\J$ to the point
$0\x\frac12\x f(\J)$ by a path $l^+_\J$ contained in
$[-1+\frac i{k+1},0]\x I\x (\frac i{k+1},\frac{i+1}{k+1}]$.
Now let us consider two close push-offs $l'_\J$ and $l''_\J$ of the
path $l_\J\cup l^+_\J$, contained in a small regular neighborhood
$L_\J$ of $l_\J\cup l^+_\J$ in the exterior of
$([-1,-1+\frac i{k+1}]\cup [0,1])\x I\x [0,\frac{i+1}{k+1}]$.
We may assume that the initial endpoints of $l'_\J$ and $l''_\J$
coincide with the two endpoints of $r_\J$.
Finally, connect the terminal endpoint of each $l'_{\I,i}$ to the terminal
endpoint of each $l''_{\I,i+1}$ by some path $c_{\I,i}$ contained in
$[-1,-1+\frac i{k+1})\x I\x(\frac i{k+1},\frac{i+1}{k+1}]$, and connect the
terminal endpoint of $l'_{\I,1}$ (resp\. $l''_{\I,n_{\I}}$) to
some point in $[-1+\frac{i-1}{k+1},-1+\frac i{k+1})\x I\x\{\frac i{k+1}\}$
(resp\. to some point in $[-1,-1+\frac{i-1}{k+1})\x I\x\{\frac i{k+1}\}$) by
some path $c_{\I,0}$ (resp\. $c_{\I,n_{\I}}$) contained in
$[-1,-1+\frac i{k+1})\x I\x [\frac i{k+1},\frac{i+1}{k+1}]$.
Define
$$l_I=\bar c_{\I,0}(\bar l'_{\I,1}r_{\I,1}l''_{\I,1}
c_{\I,1})\dots(\bar l'_{\I,n_{\I}}r_{\I,n_{\I}}
l''_{\I,n_{\I}}c_{\I,n_{\I}}).$$

In particular, this defines the required path $l$, and so it remains to
construct the balls $B_1,\dots,B_k$.
Set $$B_i=[-1,1]\x I\x [0,\frac i{k+1}]\cup
[-1,-1+\frac i{k+1}]\x I\x [\frac i{k+1},\frac{i+1}{k+1}]\cup
\bigcup_{\I\in S_{g,i}}(L_{\I}\cup D_{\I}).$$
It is easy to see that $l\i B_1\i\dots\i B_k\i S^3\but (L(mS^1)\but K)$ and
that $L(mS^1)\cap B_i$ is contained in an arc contained in $B_{i+1}$
for each $i$.
Also, without loss of generality, $o\i B_1$. \qed
\enddemo

\proclaim{Theorem 3.6} For any $k$ and $L$ the group $\G_k(L)$ is nilpotent.
\endproclaim

\remark{Remark} In the case where $k=1$ and $\pi(L)$ can be generated by two
meridians, Theorem 3.6 is due to R.~ Mikhailov, who used a thorough analysis
of the Hall basis to verify it in this case.
Since the fundamental group of the Mazur link can be generated by two
meridians, this was enough to imply that $\G_1(L)$, as defined in this section
(Theorem 3.5 was not available at that time), is not a complete invariant
of $1$-quasi-isotopy \cite{MM}.
However, since the fundamental group of the string link $W\#\rho W$ cannot be
generated by two meridians, Mikhailov's result added little with respect
to Problem 1.5.
In fact, misled by similarity between $\G_k(L)$ and groups generated by
their $(k+2)$-Engel elements (which will be discussed in
{\it ``$n$-Quasi-isotopy III''}), the first author and Mikhailov were trying
to prove that $\G_1(W\# rW)$ was not nilpotent for a considerable time;
cf\. \cite{MM}.
\endremark

\demo{Proof} We recall that $\G_k(L)=\pi(L)/\mu_k$ where
$$\mu_k=\left<[m,m^g]\mid m\in M, g\in\left<m\right>_k^{\pi(L)}\right>.$$
Using the commutator identities $[a,bc]=[a,c][a,b]^c$ and
$[a,b^{-1}]=[a,b]^{-b^{-1}}$ together with the fact that $\mu_k$ is
normal in $\pi(L)$, we see that $\mu_k$ contains any commutator of type
$[m,m^{\iota_1 g_1}\dots m^{\iota_r g_r}]$ where $m\in M$, $r$ is a positive
integer, $\iota_i=\pm 1$ and $g_i\in\left<m\right>_k^{\pi(L)}$.
Hence $\mu_k$ is the subgroup generated by
$\{[m,g]\mid m\in M, g\in\left<m\right>_{k+1}^{\pi(L)}\}$.
Analogously any commutator $[m^n,g]$, where $m\in M$, $n\in\Z$ and
$g\in\left<m\right>_{k+1}^{\pi(L)}$, can be expressed in the elements of
$\mu_k$, thus $$\mu_k=\left<\bigcup_{m\in M}[\left<m\right>,
\left<m\right>_{k+1}^{\pi(L)}]\right>.\tag{3.1}$$
Then in the quotient over $\mu_k$, each element $\bar m$ of the image $\bar M$
of $M$ under the quotient map is contained in the center of the subgroup
$\left<\bar m\right>_{k+1}^{\G_k}$.
In particular, the cyclic subgroup $\left<\bar m\right>$ is normal in
$\left<\bar m\right>_{k+1}^{\G_k}$, and thus subnormal%
\footnote{A subgroup $H$ of a group $G$ is said to be {\it subnormal} in $G$
if there exists a finite chain of subgroups $H=H_0\i H_1\i\dots\i H_d=G$ such
that each $H_i$ is normal in $H_{i+1}$.
The minimal such $d$ is called the {\it defect} of $H$.}
(of defect $\le k+2$) in $\G_k$, namely, we have that
$\left<\bar m\right>\normal\left<m\right>_{k+1}^{\G_k}\normal
\left<m\right>_k^{\G_k}\normal\dots\normal\left<m\right>_0^{\G_k}=\G_k$.

So our group $\G_k(L)$ is generated by a finite number of its subnormal cyclic
subgroups, and it remains to apply the following group-theoretic fact. \qed
\enddemo

\proclaim{Theorem (Baer)} \cite{Ba}, \cite{LS; Theorem 1.6.2} A group generated
by finitely many of its subnormal finitely generated nilpotent subgroups is
nilpotent.
\endproclaim

For convenience of the reader, we sketch a quick proof of the Baer Theorem
modulo the following textbook level result.

\proclaim{Theorem (Hirsch--Plotkin)} \cite{Rob} Let $H_1$ and $H_2$ be normal
locally nilpotent%
\footnote{A group is called {\it locally nilpotent} if all of its finitely
generated subgroups are nilpotent.
The same statement with ``locally nilpotent'' replaced by ``nilpotent''
(Fitting's Theorem \cite{Rob}) would not fit in our argument since a group
generated by its two subnormal nilpotent subgroups may be non-nilpotent
\cite{LS; p\. 22}, in contrast to the Baer Theorem.}
subgroups of a group $G$.
Then $H_1 H_2$ is a normal locally nilpotent subgroup of $G$.
\endproclaim

Clearly, the union of ascending chain of locally nilpotent subgroups is locally
nilpotent, so using the Zorn Lemma (for our purpose, the countable case
suffices) one obtains that in any group $G$ there is a unique normal
locally nilpotent subgroup (called the Hirsch--Plotkin radical of $G$)
containing all normal locally nilpotent subgroups of $G$.
Using an induction on the defect and applying, at each step, the
Hirsch--Plotkin Theorem (which immediately generalizes to the case of any
(ordinal) number of subgroups), one proves (cf\. \cite{Rob}) that the
Hirsch--Plotkin radical contains all the subnormal locally nilpotent subgroups.
In particular, we obtain that any group generated by its subnormal locally
nilpotent subgroups is locally nilpotent.

\remark{Remark} In connection with Milnor's original definition of
$\G(L)=\G_0(L)$ \cite{Mi}, we notice that $\mu_k$ is the subgroup generated by
the commutator subgroups of the subgroups $\left<m\right>_{k+1}^{\pi(L)}$
where $m$ runs over all meridians.
Indeed, by the definition, $\left<m\right>_{k+1}^{\pi(L)}$ is generated by the
conjugates of $m$ under the action of $\left<m\right>_k^{\pi(L)}$, and if $m^g$
and $m^h$ are two of them, $[m^g,m^h]=[m,m^{hg^{-1}}]^g\in\mu_k$.
Conversely, formula (3.1) implies that $\mu_k$ is contained in the subgroup
generated by these commutators.
\endremark

\proclaim{Corollary 3.7} Every quotient of the fundamental group of a link
(or string link), functorially invariant under strong $k$-quasi-isotopy for
some $k$, is nilpotent.
\endproclaim

The case of links follows from Theorems 3.5 and 3.6; the proof in the case of
string links is analogous.

\head Appendix. A proof that $\tl\beta\equiv\bar\mu(1122)\pmod\lk$ \endhead

We give a direct proof  that the Milnor invariant $\bar\mu(1122)\in\Z/\lk$
satisfies the crossing change formula (2.2).
It is easy to check that $\bar\mu(1122)$ vanishes on the links $T_n$ from
\cite{KL1}, hence coincides with the residue class of $\tl\beta\bmod\lk$
(which also follows from \cite{Co3}).

\medskip
We recall the definition of $\bar\mu$-invariants following Levine's approach
\cite{Le}.
Given a link $L'=L_1\sqcup\dots\sqcup L_n\i S^3$, we assume inductively
that the $\bar\mu$-invariants of $L=L_1\sqcup\dots\sqcup L_{n-1}$ are already
defined.
To define those of $L'$, we only need to know the homotopy class of $L_n$ in
$\pi(L)=\pi_1(S^3\but L)$ (up to conjugation).
Furthermore, we proceed to the quotient $\G(L)=\pi(L)/\mu_0$ (see \S3).
A presentation for $\G(L)$ is (cf\. \cite{Le})
$$\left<m_1,m_2,\dots,m_{n-1}\mid [m_i,l_i]=1, R_j=1\right>,$$
where $m_1,\dots,m_{n-1}$ is a collection of meridians, one to each component,
$l_1,\dots,l_{n-1}$ are the corresponding longitudes (regarded as words in
$m_i$'s), and $R_j$, $j=1,2,\dots$ are all commutators in $m_i$'s with repeats,
i.e\. with some $m_i$ occurring at least twice.

Now any element $\alpha\in\G(T)$, where $T$ is the trivial link and therefore
$\pi(T)$ is the free group $F_{n-1}$ with free generators $m_1,\dots,m_{n-1}$,
has unique decomposition into product of powers of {\it basic commutators
without repeats} (see \cite{Le} for definition and proof; for a general
treatment of basic commutators see \cite{Hall}), which is finite since $\G(L)$
is nilpotent by \cite{Mi1}.
For our case $n=4$ the decomposition looks like $$\alpha=x^{e_1}y^{e_2}z^{e_3}
[y,x]^{e_4}[z,y]^{e_5}[z,x]^{e_6}[[y,x],z]^{e_7}[[z,x],y]^{e_8},\tag{A.1}$$
where $m_1,m_2,m_3$ are substituted for convenience by $x,y,z$ respectively.
The integers $e_i$ depend on the choice of basis in $F_{n-1}$, and hence
should be considered only up to certain transformations (all of them are
generated by conjugations $m_i\mapsto (m_i)^g$, see \cite{Le} for details).
Next, a choice of meridians yields a homomorphism $\pi(T)=F_{n-1}\to\pi(L)$
which descends to $\phi\:\G(T)\to\G(L)$, so for arbitrary $L$ and
$\alpha\in\G(L)$ we are led to the representation (A.1).
However, in general $l_i$ are not trivial and so the relations $[m_i,l_i]=1$
arising from the peripheral tori are not vacuous, hence the representation
(A.1) is not unique with respect to a fixed basis of $m_i$.
The commutator numbers $e_i=(e_i)_j$ for $\alpha=\alpha_j=[m_j,l_j]$ give the
precise indeterminacy in this case.
Finally, if we set $\alpha$ to be the class of $L_n$ in $\G(L)$, the obtained
commutator numbers $e_i$ up to the prescribed indeterminacy are, together with
the invariants of $L$, the complete set of homotopy invariants of $L'$
\cite{Le}, including all Milnor's invariants with their original
indeterminacy.
In particular, in the formula (A.1) above, $e_1,e_2,e_3$ are the linking
numbers of $L_4$ with $L_1$, $L_2$, $L_3$ respectively, the residue classes
of $e_4,e_5,e_6$ are the triple $\bar\mu$-invariants corresponding to
$(L_i,L_j,L_4)$, where $(i,j)$ run over the two-element subsets of $1,2,3$,
and the residue classes of $e_7,e_8$ are the two independent quadruple
invariants $\bar\mu(1234)$, $\bar\mu(1324)$, respectively.

Now let $\bar L_\pm=K_1^\pm\sqcup K_3^\pm$ be two-component links of linking
number $l$, related by a single positive crossing change so that the formula
(2.2) applies.
Let $K_2^\pm$ and $K_4^\pm$ be zero push-offs of $K_1^\pm$ and $K_3^\pm$,
respectively, and let us consider the links $L_\pm:=\bar L_\pm\sqcup K_2^\pm$
and $L_\pm':=L_\pm\sqcup K_4^\pm$.

Let us see the difference between the commutator numbers for the elements
$\alpha_-=[K_4^-]\in\G(L_-)$ and $\alpha_+=[K_4^+]\in\G(L_+)$.
Since $L_-$ (Fig\. ~5a) and $L_+$ (Fig\. ~5b) are link homotopic (by the
obvious homotopy with a single self-intersection on $K_3^\pm$), the two groups
$\G(L_-)$, $\G(L_+)$ are isomorphic by an isomorphism $h$ taking
meridians onto meridians and longitudes algebraically somewhat close to
longitudes \cite{Mi1; Theorem 2}, \cite{Le; Theorem 4} or Theorem 3.2.
In our case the only longitude one has to worry about is $\alpha_-$, and its
behavior can be determined.
Let $\tau_1,\tau_2\in\G(L_+)$ denote the classes of the two lobes' longitudes,
joined to the basepoint as shown in Fig\. 5b, so that the linking number of
each $\tau_i$ with $K_3^+$ is zero.
By Lemma 3.3, $h(\alpha_-)=\xi\alpha_+$, where
$\xi=[\tau_1,[\tau_1^{-1},z^{-1}]]$.
We will systematically use the fact that the fourth central term
$\gamma_4\G(L_+)=1$ since $L_+$ has three components \cite{Mi1} allowing us
to use simplified commutator identities, e.g\. it follows that
$\xi=[[z,\tau_1],\tau_1]$.
In fact, as we noted above, $\xi$ is not well defined as a word in meridians,
and can be rewritten using the relations $[m_i,l_i]=1$ arbitrarily.
In particular, since $\alpha_+=z\tau_1 z^{-1}\tau_2$, the relation
$[z^{-1},\alpha_+]=1$, or, equivalently, $[\tau_2^{-1},z]=[z^{-1},\tau_1]$,
implies that $\xi$ can be rewritten e.g\. as $\xi'=[[\tau_2,z],\tau_1]$.

\fig 5

Now suppose that the linking numbers of the lobes with the first component of
$\bar L_+$ are $n$ and $l-n$.
Then the representations (A.1) of $\tau_i$'s are as follows:
$$\tau_1=x^n y^n\dots;\qquad \tau_2=x^{l-n}y^{l-n}\dots$$
The higher terms are omitted since they will not affect the result, as all
commutators with repeats or of weight four are trivial.
Making the substitutions for $\tau_i$'s, we obtain (using the Hall--Witt
identity $[[x,y],z][[y,z],x][[z,x],y]=1\bmod\gamma_4$)
$$\eightpoint \xi'=[[x^{l-n}y^{l-n},z],x^n y^n]=
[[y,z],x]^{n(l-n)}[[x,z],y]^{n(l-n)}=[[y,x],z]^{n(l-n)}[[z,x],y]]^{-2n(l-n)}$$
which is not surprising, since $\bar\mu(1212)=-2\bar\mu(1122)$ \cite{Mi}, and
proves the crossing change formula (2.2) for $\bar\mu(1122)$.
The analogous calculation for $\xi$:
$$\xi=[[z,x^n y^n],x^n y^n]=[[y,x],z]^{-n^2}[[z,x],y]^{2n^2}$$
shows that the integer $e_7$ is not well-defined even for the two-component
link.

\head Acknowledgements \endhead

The first author is indebted to the referee of {\it ``$n$-Quasi-isotopy III''}
for pointing out an error in this paper, and would like to thank M. Cencelj,
A. N. Dranishnikov, U. Kaiser, J. Male\v{s}i\v{c}, V. M. Nezhinskij,
R. Sadykov, E.~ V.~ Shchepin, A. Skopenkov, A. Turull and especially
R. Mikhailov and P. Akhmetiev for stimulating conversations.
We are grateful to the referee of this paper for a thoughtful reading and
to the editors for kindly allowing us to include inappropriately numerous and
late corrections and elaborations.
Some results of this paper and {\it ``$n$-Quasi-isotopy II''} were announced
in \cite{MM}.
The authors are partially supported by the Ministry of Science and Technology
of the Republic of Slovenia grant No\. J1-0885-0101-98 and the Russian
Foundation for Basic Research grants 99-01-00009, 02-01-00014 and 05-01-00993.

\enddocument
\end